\documentclass[a4paper]{amsart}
\usepackage{amsthm}
\usepackage{amssymb}
\usepackage{amsmath}
\usepackage{verbatim,enumerate,graphicx,wrapfig}

\renewcommand{\thefootnote}{\fnsymbol{footnote}}

\newtheorem{theorem}{Theorem}[section]
\newtheorem{lemma}[theorem]{Lemma}

\newtheorem{proposition}[theorem]{Proposition}
\newtheorem{remark}[theorem]{Remark}

\newtheorem{example}[theorem]{Example}

\newtheorem*{example*}{Example}

\newtheorem*{remark*}{Remark}

\newtheorem{corollary}[theorem]{Corollary}
\newtheorem*{corollary*}{Corollary}

\newtheorem{definition}[theorem]{Definition}

\newtheorem*{definition*}{Definition}

\newtheorem*{notation*}{Notation}

\newtheorem{notation}[theorem]{Notation}

\numberwithin{equation}{section}

\gdef\myletter{}

\let\savetheequation\theequation
\def\theequation{\savetheequation\myletter}

\def\bt{\mathbf{t}}
\def\bp{\mathbf{p}}

\newcommand{\CC}{{\mathbb C}}

\newcommand{\RR}{{\mathbb R}}
\newcommand{\ZZ}{{\mathbb Z}}

\newcommand{\PP}{{\mathbb P}}

\newcommand{\NN}{{\mathbb N}}

\newcommand{\be}{\mathbf{e}}

\newcommand{\calP}{\mathcal{P}}
\newcommand{\calM}{\mathcal{M}}

\newcommand{\calE}{\mathcal{E}}
\newcommand{\calN}{\mathcal{N}}
\newcommand{\calU}{\mathcal{U}}
\newcommand{\lt}{\textsc{lt}}
\newcommand{\lc}{\textsc{lc}}

\newcommand{\ttt}{\textsc{tt}}
\newcommand{\tc}{\textsc{tc}}
\newcommand{\tm}{\textsc{tm}}

\newcommand{\longmapsfrom}{\mathrel{\reflectbox{\ensuremath{\longmapsto}}}}

\newcommand{\vol}{\mbox{vol}}

\def \hat{\widehat}

\def \b0{{\bf 0}}
\def\bV{{\bf V}}
\def\bI{{\bf I}}

\def \span{\mathrm{span}}

\def\calB{\mathcal{B}}

\def\calO{\mathcal{O}}

\long\def\symbolfootnote[#1]#2{\begingroup%
\def\thefootnote{\fnsymbol{footnote}}\footnote[#1]{#2}\endgroup}

\begin{document}

\title{Okounkov bodies and transfinite diameter}
\author{Sione Ma`u}
\address{Department of Mathematics,
University of Auckland,
Auckland, NZ}
\email{s.mau@auckland.ac.nz}
\begin{abstract}We present an explicit calculation of an Okounkov body associated to an algebraic variety.  This is used to derive a formula for transfinite diameter on the variety.  We relate this formula to a recent result of Nystr\"om.
\end{abstract}
\keywords{Chebyshev transform, transfinite diameter, algebraic variety, Okounkov body, polynomial}
\subjclass[2010]{32U20; 14Q15}

\maketitle

\section{Introduction}

This paper investigates Okounkov bodies and their relation to transfinite diameter.   
 Lazarsfeld and Mustata \cite{lazarsfeldmustata:convex}, as well as Kaveh and Khovanskii \cite{kavehkhovanskii:newton},  introduced Okounkov bodies into algebraic geometry as an important tool for the asymptotic study of linear series.\footnote{Think of these as classes of polynomials.}  At about the same time, Berman and Boucksom \cite{bermanboucksom:growth} used pluripotential theory to study the asymptotic properties of powers of big line bundles.  These powers form a natural linear series.  Therefore Okounkov bodies ought to be related to pluripotential theory.  Such a connection was made in a theorem of Nystr\"{o}m \cite{nystrom:transforming}.

Nystr\"{o}m's result on Okounkov bodies is closely related to a classical theorem of Zaharjuta \cite{zaharjuta:transfinite}.  The latter (reproduced in this paper as Theorem \ref{thm:z}) gives an important capacity in pluripotential theory---the transfinite diameter of a compact set $K\subset\CC^n$---as a real integral over an $(n-1)$-dimensional simplex in $\RR^n$ involving so-called \emph{directional Chebyshev constants}, which are quantities defined in terms of polynomials on $K$. 
 Nystr\"{o}m gives a similar looking integral formula that relates the Monge-Amp\`ere energy of Hermitian metrics on a line bundle $L$ over a compact complex manifold, to an integral over the Okounkov body associated to $L$ of the Chebyshev transforms of these metrics.  (We will explain the  terms later.) 

To relate the two results, classical objects in pluripotential theory need to be transferred into the modern theory on complex manifolds:
\begin{eqnarray}
\hbox{polynomials} &\hookrightarrow& \hbox{sections of line bundles}  \label{eqn:ps} \\
    \hbox{compact set} &\hookrightarrow&  \hbox{Hermitian metric}   \label{eqn:ch}  \\
 \hbox{simplex} &\hookrightarrow&  \hbox{Okounkov body}    \label{eqn:so}
 \end{eqnarray}
where $A\hookrightarrow B$ means that all objects in $A$ can be modelled (more or less) by objects in $B$.  

The main aim of this paper is to prove a version of Nystr\"{o}m's result on algebraic subvarieties of $\CC^n$.  A subvariety of $\CC^n$ is an intermediate setting which is both a natural extension of the classical theory in $\CC^n$ as well as a concrete illustration of the theory on complex manifolds.  It provides a natural bridge from the classical to the modern point of view.

    Zaharjuta's methods may be naturally adapted to this setting, with some additional tools from computational algebraic geometry and {weighted} pluripotential theory.  The methods used in this paper are also similar to \cite{coxmau:transfinite}, especially the use of computational algebraic geometry to carry out computations on a variety.  They do not use much pluripotential theory: no plurisubharmonic functions or Monge-Amp\`ere integrals are required, only polynomials.

We begin with some background material and describe the relationships (\ref{eqn:ps}),  (\ref{eqn:ch}), and (\ref{eqn:so}) given above.  Section \ref{sec:prelim} deals with the first two.  First, (\ref{eqn:ps}) is elementary complex geometry; one may skip this part if one is familiar with relating polynomials on a  variety in $\CC^n$ to powers of the line bundle $\calO(1)$ over the corresponding variety in $\PP^n$. Next, to describe (\ref{eqn:ch}) we use the notion of a \emph{(weakly) admissible  weight}.  The identification of Hermitian metrics and weights is reasonably familiar from the application of  pluripotential theory to complex geometry (see e.g. \cite{demailly:singular}, \cite{guedjzeriahi:intrinsic}).  Capacities associated to sets can be defined in terms of weights supported on these sets.\footnote{If the weight is identically 1 we recover the classical (unweighted) theory.} 

In section \ref{sec:ok} we define the Okounkov body and study some of its properties.  This definition depends on a choice of coordinates, and we use {Noether normalization from computational algebraic geometry to choose coordinates} which have good computational properties.   It is easily seen by definition that the Okounkov for $\CC^n$ is the region $S$ bounded by the coordinate hyperplanes and the standard simplex in $\RR^n$. 
 (This provides, more or less, the connection (\ref{eqn:so}).)   We then present a fairly explicit algorithm for constructing an Okounkov body associated to a variety in $\CC^n$.  Although we do not give a rigorous proof of the method in general, we use it to compute the Okounkov body associated to a complexified unit sphere in $\CC^3$.  

In section \ref{sec:dir}, we define directional Chebyshev constants, the Chebyshev transform, and a notion of transfinite diameter.  We then prove our main theorem (Theorem \ref{thm:Ok}) that gives transfinite diameter on the variety as an integral over the Okounkov body.  The relation with Nystr\"om's result is seen by translating things into the language of complex geometry.

In section \ref{sec:comparison} we make the explicit connection to Zaharjuta's classical theorem, as well as a  homogeneous version of Jedrzejowski \cite{jedrzejowski:homogeneous}.  This involves making a projective change of coordinates.  (It is good to have the complex geometric point of view here.)

Finally, in section \ref{sec:further} we investigate further properties of Chebyshev constants on the sphere in $\CC^3$.  In particular, we look at directional Chebyshev constants associated to so-called \emph{locally circled sets} (Proposition 6.7).  The notion of a locally circled set is adapted from the classical notion of a circled set in $\CC^n$ (cf. \cite{bloom:weighted}, \cite{bloomlev:weighted}). 

\section{Preliminaries} \label{sec:prelim}
\subsection{Varieties in $\CC^n$} Let $\CC[z]=\CC[z_1,\ldots,z_n]$ denote the ring of polynomials in $n$ variables.  Recall that an algebraic variety $V$ in $\CC^n$ ($n>1$ is an integer), is the solution to a finite collection of polynomial equations
$$
V=\{a\in\CC^n: P_1(a)=\cdots=P_m(a)=0, \ P_j\in\CC[z] \forall\, j\};
$$
\begin{notation}\rm
Given an algebraic variety $V\subseteq\CC^n$, define 
$$\bI(V):= \{p\in\CC[z]:\ p(a)= 0 \hbox{ for all } a\in V\}.$$
It is easy to see that this is an ideal.  Also, given an ideal $I\subseteq\CC[z]$, define
$$\bV(I):=\{a\in\CC^n:\ p(a)=0 \hbox{ for all } p\in\CC[z]\}.$$
\end{notation}
\begin{theorem} \rm
\begin{enumerate}
\item (Hilbert basis theorem) Any ideal $I$ is finitely generated; consequently, $\bV(I)$ is always an algebraic variety.
\item (Nullstellensatz) For any ideal $I\subseteq\CC[z]$,  we have $I\subseteq\bI(\bV(I))$, and if  the property \begin{equation}\label{eqn:radical} p^m\in I \hbox{ for some } m\in\NN \Longrightarrow p\in I,\end{equation} holds, then $I=\bI(\bV(I))$.  
\item For any algebraic variety $V$,  $\bI(V)$ satisfies (\ref{eqn:radical}) and $\bV(\bI(V))=V$.  \end{enumerate}
\end{theorem}
Suppose $V=\bV(I)$ where $I=\langle P_1,\ldots,P_m\rangle$ is the ideal generated by the polynomials $P_j$.  If $I$ satisfies (\ref{eqn:radical}) then the above theorem implies that restricting the evaluation of polynomials to points of $V$ is equivalent to working with elements of the factor ring $\CC[z]/I$ via the correspondence 
$$p=q \hbox{ on } V\iff p-q\in I. 
$$

\subsection{Projective space} Algebraic subvarieties of $\CC^n$ can be put into the complex geometric setting using projective space. Consider $\CC^n\subset\PP^n$ via the usual embedding $$z=(z_1,\ldots,z_n)\hookrightarrow [1:z_1:\cdots:z_n]=[1:z],$$
where we use homogeneous coordinates on the right-hand side: $\PP^n=\CC^{n+1}/\sim$ with the equivalence $(z_0,\cdots,z_n)\sim(w_0,\cdots,w_n)$ if there is a $\lambda\in\CC$ such that $\lambda z_i=w_i$ for each $i$; we write $[z_0:z]=[z_0:\cdots:z_n]=[w_0:\cdots:w_n]=[w_0:w]$.  We have $\PP^n= \CC^n\cup H_{\infty}$ where $H_{\infty} = \{[0:z]\in\PP^n: z\in\CC^n\}$ is the hyperplane at infinity.  

 The standard affine charts of $\PP^n$ as a complex manifold will be denoted by $\calU_j$, $j=0,\ldots,n$. These are given by $\calU_0=\CC^n$ with the standard embedding described above, and for $j>0$, 
$
\calU_j=\{[Z_0:Z_1:\cdots:Z_N]\in\PP^n: Z_j\neq 0\} 
$
with the map 
\begin{equation}\label{eqn:21}
\CC^n\ni(w_1,\ldots,\hat w_j,\ldots,  w_n)\hookrightarrow[w_1:\cdots:w_{j-1}:1:w_{j+1}:\cdots:w_n]\in\calU_j
\end{equation}
 giving local coordinates on $\calU_j$ (here $\hat w_j$ means that there is no $w_j$ coordinate).  Going the other way is \emph{dehomogenization}: $$\calU_j\ni[Z_0:\cdots:Z_n]\mapsto (Z_0/Z_j,\ldots,Z_n/Z_j)\in\CC^n.$$  We also have the change of coordinates on the overlap $\calU_j\cap\calU_k$:
$$
w_{j}=1/v_k,\quad  w_{\ell}=v_{\ell}/v_k \ \hbox{for all }\ell\neq j, 
$$
with $(v_0,\ldots,\hat v_j,\ldots,v_n)\in\calU_j$ and $(w_0,\ldots,\hat w_k,\ldots,w_n)\in\calU_k$. 

Let $V_{\PP}\subset\PP^n$ be the continuous extension of $V\subset \CC^n$ across points of $H_{\infty}$ under the above embedding (the \emph{projective closure}).  One can use homogeneous coordinates to characterize it:   
$$
V_{\PP} = \{[Z_0:\cdots:Z_n]\in\PP^n: \ p(Z_0,\ldots,Z_n)=0 \hbox{ for all } p\in \bI_h(V)  \}
$$
where $\bI_h(V)\subset\CC[Z_0,\ldots,Z_n]$ is the collection of homogeneous polynomials $p(Z)$ such that $p(1,z_1,\ldots,z_n)= 0$ whenever $(z_1,\ldots,z_n)\in V$.

\subsection{Sections of line bundles and polynomials}
Let us recall the basic notions associated to a holomorphic line bundle $L$ over a complex manifold $M$ of dimension $m$, which is essentially a union of complex lines (i.e. complex vector spaces of dimension 1 or copies of $\CC$) parametrized holomorphically by points of $M$.\footnote{In what follows all notions will, unless otherwise stated, refer to their complex versions.} Precisely, $L$ is a manifold of dimension $m+1$ with a projection $\pi:L\to M$ such that $L_a:=\pi^{-1}(a)$ is a complex line for each $a\in M$.  A \emph{(holomorphic) section of $L$} is a holomorphic map $s:M\to L$ with $(\pi\circ s)(a)=a$.

We review the details of the local product structure of $L$: any point $a\in M$ has a neighborhood $U$ for which there is a holomorphic injection 
$U\times\CC\ni(z,\zeta)\stackrel{\varphi}{\mapsto} v\in  L$
such that for each $z\in U$,  
$\pi\circ\varphi(z,\zeta)=z$ and the map $\CC\ni\zeta\mapsto\varphi(z,\zeta)\in L_a$ {is a linear isomorphism}.  The pair $(U,\varphi)$ is called a \emph{local trivialization}.  

Let $\{(U_{\alpha},\varphi_{\alpha})\}_{\alpha}$ be a collection of local trivializations that cover $L$, i.e., $\bigcup_{\alpha} U_{\alpha}= M$, so that $\bigcup_{\alpha}\left(\bigcup_{a\in U_{\alpha}} L_a\right)=L$.  Define 
$g_{\alpha\beta}:U_{\alpha}\cap U_{\beta}\to\CC$ by $g_{\alpha\beta}(z)=\eta/\zeta$, where 
$$
U_{\alpha}\times\CC \ni (z,\zeta) \stackrel{\varphi_{\alpha}}{\longmapsto} \  v \  \stackrel{\varphi_{\beta}}{\longmapsfrom} (w,\eta) \in U_{\beta}\times\CC  
$$
for some $v\in \bigcup_{a\in U_{\alpha}\cap U_{\beta}}  L_a$.  Clearly for any $\alpha,\beta,\gamma$ and $z\in U_{\alpha}\cap U_{\beta}\cap U_{\gamma}$, 
$$
g_{\alpha\beta}(z)g_{\beta\alpha}(z) = g_{\alpha\beta}(z)g_{\beta\gamma}(z)g_{\gamma\alpha}(z)= 1
$$
(the \emph{cocycle condition}.)  


Given a section $s$ of $L$, there is an associated collection $\{s_{\alpha}\}_{\alpha}$ of (local) functions  $s_{\alpha}:U_{\alpha}\to\CC$ such that
$$
s_{\beta}(z) = g_{\alpha\beta}(z)s_{\alpha}(z), \quad \hbox{for all } z\in U_{\alpha}\cap U_{\beta}.
$$
It is straightforward to verify that a collection $\{s_{\alpha} \}_{\alpha}$ of functions satisfying the above conditions characterizes a section $s$ (since $L$ is determined by $\{g_{\alpha\beta}\}$).

\medskip

We now specialize to our context.  Define $\calO(1)$ as the collection of pairs 
\begin{equation}\label{eqn:22}
\calO(1) = \{([Z_0:\cdots:Z_n],a_0Z_0+\cdots+a_nZ_n): [Z_0:\cdots:Z_n]\in\PP^n, \ a_j\in\CC    \}.
\end{equation}
The line bundle structure of $\calO(1)$ comes from function evaluation. Let us see how this works by fixing $Z\in\PP^n$ and computing $L_Z$ explicitly. First, pick $b(Z)$ such that $b(Z)\neq 0$; this is true (or not) independently of the homogeneous coordinates used to compute $b(Z)$.  We claim that
\begin{equation}\label{eqn:23a}
L_Z = \{(Z,\lambda b(Z)): \lambda\in\CC\}.
\end{equation}
For any $a(Z)=a_0Z_0+\cdots+a_nZ_n$, define $\lambda\in\CC$ by $\lambda:=\frac{a(Z)}{b(Z)}$; note that this computation of $\lambda$ is independent of homogeneous coordinates.  Rewrite this as  $a(Z)=\lambda b(Z)$, and substitute into (\ref{eqn:22}) to get (\ref{eqn:23a}).  This also verifies that $L_Z$ is indeed a complex line.

The sections of $\calO(1)$ can be immediately read off from (\ref{eqn:22}) as the objects $a(Z)$, identified with linear homogeneous polynomials in $n+1$ variables.  They form a space of dimension $n+1$, usually denoted by $H^0(\PP^n,\calO(1))$.  Now  $\CC[z_1,\ldots,z_n]_{\leq 1}$, the polynomials of degree at most 1 in $n$ variables, can be mapped into $H^0(\PP^n,\calO(1))$ by homogenizing coordinates, 
$$
a_0+a_1z_1+\cdots+a_nz_n = a(z) \longmapsto  (Z,a(Z)) \in H^0(\PP^n,\calO(1))
$$
where $Z=[Z_0:\cdots:Z_n]=[1:z_1:\cdots:z_n]$  (see (\ref{eqn:21})).  

For fixed $z\in\CC^n$ (and associated $Z=[1:z]\in\calU_0\subset\PP^n$), it is an exercise to show that this identifies $a(z)$ as the local function on $\calU_0$ of the section given by $a(Z)$ under the local trivialization 
$$\CC^n\times\CC\ni (z,\lambda) \mapsto (Z,\lambda Z_0)\in\calO(1),$$
where the right-hand side is as in (\ref{eqn:23a}), with $\lambda=a(Z)/Z_0$.  

For other values of $j$, a similar formula holds; e.g. when $j=1$, consider $$a(w)=a_0w_0+a_1+a_2w_2+\cdots+a_nw_n\in\CC[w_0,w_2,\ldots,w_n]_{\leq 1}.$$ Form $a(Z)$ with $Z=[w_0:1:w_2:\cdots:w_n]$; then $a(w)$ corresponds to $a(Z)$ under $(w,\lambda)\mapsto (Z,\lambda Z_1)$. 

We can also calculate the maps $g_{jk}$ on $\calU_j\cap\calU_k$; let us do the case $j=0,k=1$.  Suppose $Z\in\calU_0\cap\calU_1$, with coordinates $z$ on $\calU_0$ and $w$ on $\calU_1$ given by  $$[1:z_1:\cdots:z_n]=Z=[w_0:1:w_2:\cdots:w_n],$$ so that $w_0=1/z_1$ and $w_j=z_j/z_1$ for $j\neq 0,1$.  Then with
$$
\zeta = a_0+a_1z_1+\cdots z_nz_n \longmapsto   (Z,a(Z)) \longmapsfrom a_0w_0+a_1+a_2w_2+\cdots+a_nw_n = \eta,
$$
we have 
\begin{eqnarray*}
g_{01}(Z)=\frac{\eta}{\zeta} &=& \frac{a_0w_0+a_1+a_2w_2+\cdots+a_nw_n}{a_0+a_1z_1+\cdots z_nz_n} \\ &=& \frac{a_0/z_1+a_1+a_2z_2/z_1+\cdots+a_nz_n/z_1}{a_0+a_1z_1+\cdots z_nz_n} = \frac{1}{z_1}.
\end{eqnarray*}
So $g_{01}(Z)=1/z_1 = w_0$. 

\medskip

For any positive integer $k$, define $\calO(k)$ to be the line bundle over $\PP^n$ given by 
$$
\calO(k) = \Bigl\{(Z,P(Z)): Z\in\PP^n, P(Z)=\sum_{|\alpha|=k}a_{\alpha}Z^{\alpha},\ a_{\alpha}\in\CC \Bigr\}
$$
where we use the standard multi-index notation $Z^{\alpha}=Z_0^{\alpha_0}\cdots Z_n^{\alpha_n}$.  Similar calculations as above yield the following:
\begin{enumerate}
\item Given $Z\in\PP^n$, the line $L_Z$ is generated by any $P(Z)$ which evaluates to a nonzero complex number;
\item When $Z=[1:z]$, the map $$\CC[z]_{\leq k}\ni p(z) \to (Z,P(Z))\in H^0(\PP^n,\calO(k))$$ identifies $\CC[z]_{\leq k}$ with $H^0(\PP^n,\calO(k))$ under $(z,\lambda)\mapsto(Z,\lambda Z_0^k)$;
\item When $Z=[1:z]=[w_0:1:w_2:\cdots:w_n]$, we have $\displaystyle g_{01}(Z) = \frac{1}{z_1^k} = w_0^k$.  
\end{enumerate}

\begin{remark}\rm 
The last item shows that $\calO(k)$ can be identified with $\calO(1)^{\otimes k}$, the $k$-th tensor power of $\calO(1)$, where we take the tensor power fiberwise (over each line $L_a$). It is an exercise  to show that the same transition functions are obtained, and hence these line bundles have the same structure. 
\end{remark}

The constructions of $\calO(k)$ on $\PP^n$ give line bundles on holomorphic submanifolds of $\PP^n$, by restriction; in particular, when $V\subset\CC^n$ is a smooth algebraic subvariety with extension $V_{\PP}\subset\PP^n$.  The restriction of $\calO(k)$ to points over $V_{\PP}$ corresponds to the restriction to $V$ of the polynomials $\CC[z]_{\leq k}=\{p\in\CC[z]: \deg(p)\leq k\}$.  

\subsection{Hermitian metrics and weights}
Recall that a Hermitian inner product on a complex vector space $V$ is a map $\langle\cdot,\cdot\rangle:\CC\times\CC\to\CC$ for which $z\mapsto\langle z,w\rangle$ is linear for each fixed $w$, and $w\mapsto\langle z,w\rangle$ is conjugate-linear for each fixed $z$, and $\langle z,w\rangle=\overline{\langle w,z\rangle}$.  

A Hermitian metric on a line bundle $L$ over $M$ is a family of Hermitian inner products $\langle\cdot,\cdot\rangle_a$ on the fibers $L_a$ varying continuously in $a$: for any sections $s,t$ the map  $a\mapsto\langle s(a),t(a)\rangle_a$ is continuous.  (In what follows we will suppress the dependence of the inner product on $a\in M$.)

\medskip

Let us specialize to the line bundles $\calO(k)$ over $\PP^n$, with sections given by homogeneous polynomials as above. Later, we will restrict to subvarieties.

\begin{example}  \label{ex:23} \rm
Consider the metric on $\calO(k)$ given by
$$\langle p(Z),q(Z)\rangle = \frac{p(Z)\overline{q(Z)}}{|Z|^{2k}},$$
where $|Z|^2=|Z_0|^2+\cdots+|Z_N|^2$, and we evaluate the right-hand side in homogeneous coordinates.  (Note that the value obtained is independent of homogeneous coordinates.)  In local coordinates on $\calU_0$, one can write this as%
\begin{equation}  \label{eqn:ex15}
\langle p(z),q(z)\rangle = \frac{p(z)\overline{q(z)}}{(1+|z|^{2})^k}.
\end{equation}
\end{example}

\begin{example}\label{ex:24} \rm 
Suppose $W:\CC^{n+1}\to\CC$ is a continuous function with the property that $|W(\lambda Z)| = |\lambda|^k|W(Z)|$. Then the formula
\begin{equation}\label{eqn:23}
\langle p(Z),q(Z)\rangle_W := \frac{p(Z)}{W(Z)} \overline{\left(\frac{q(Z)}{W(Z)}\right)}
\end{equation}
defines a Hermitian metric on $\calO(k)$.  

Replacing $W$ by its absolute value $|W|$ makes no difference to the right-hand side, but if $W$ is holomorphic in some region it might be useful to leave this structure intact.
\end{example}

As before, sections of $\calO(k)$ may be identified with polynomials in $\CC[z]_{\leq k}$ by transforming to affine coordinates on $\calU_0$; we will see below that the Hermitian metric may be identified with a weight on $\CC^n$.  We now define what this is.

\begin{definition}\rm
Let $K\subset\CC^n$ be a set.  An \emph{weight function on $K$} is a function $w:K\to\CC$ for which
\begin{enumerate}\item the absolute value $z\mapsto|w(z)|$ is lower semicontinuous; and
\item there is a non-negative real number $r$ such that $|w(z)|$ decays like $o(|z|^{-r})$ as $|z|\to\infty$.  Let us denote by $r(w)$ the inf over all such $r$.
\end{enumerate}
If $r(w)<1$ then $w$ is said to be an \emph{admissible weight function}.  If $r(w)\leq 1$, then $w$ is \emph{weakly admissible}. 
\end{definition}
Clearly  $r(w^t)=t r(w)$ for any positive integer $t$.

\begin{example}\rm 
Consider a polynomial $p\in\CC[z]$ of degree $d\in\NN$.  For any $\epsilon>0$, the function $w={1}/{p}$ is a weight on any set   $K\subseteq(\CC^n\setminus\{z:|p(z)|>\epsilon\})$. If $K$ is bounded then $r(w)=0$, otherwise  $r(w)=d$, and the weight given by $|p|^{-1/d}$ is weakly admissible.  \end{example}

\medskip

An admissible weight function on $K\subset\CC^n$ is used to evaluate polynomials.
\begin{definition} \label{def:27}  \rm
Let $w:K\to\CC$ be an admissible weight on $K$.  
For any $p\in\CC[z]$ we define the \emph{weighted polynomial evaluation}
$$
p(z)_w:=w(z)p(z),\quad   p(z)_{w,k}:=  w(z)^kp(z) \quad (z\in K,\ k\in\NN),
$$
which we extend by zero: $p(z)_{w}=p(z)_{w,k}=0$ if $z\not\in K$.  This also yields the \emph{weighted sup norms} 
$$
\|p\|_{K,w,k}:=  \|w^kp\|_K = \sup_{z\in K} |w(z)^kp(z)| \quad (k\in\NN).
$$
\end{definition}

Let us relate the Hermitian metric on $\calO(k)$ given by Example \ref{ex:24}  to a weight on $\CC^n$. Writing equation (\ref{eqn:23}) in terms of $z$ coordinates, where $[1:z]=Z\in\calU_0$, we have 
$$
\langle p(Z),q(Z)\rangle_W := 
\frac{p(z)}{W([1:z])} \overline{\left(\frac{q(z)}{W([1:z])}\right)} = p(z)_{w,k}\overline{q(z)_{w,k}}
$$
where $w(z):={|W([1:z])}|^{-1/k}$ defines the weight function.  If we compare the above to (\ref{eqn:ex15}) in  Example \ref{ex:23}, 
$$
\left|\frac{\langle p(Z),q(Z)\rangle_W}{\langle p(Z),q(Z)\rangle} \right| = \frac{|w(z)|^{2k}}{(1+|z|^2)^k}.
$$
Using the above formula we can deduce that $w$ must be continuous and weakly admissible.  For fixed $Z$, we can choose $p,q$ for which $\langle p(Z),q(Z)\rangle\neq 0$ in a neighbourhood of $Z$.  Then the left-hand side is continuous.  The right-hand side then shows that this quantity is independent of $p$ and $q$ on $\CC^n=\PP^n\setminus H_{\infty}$, and hence on $\PP^n$ (extending by continuity).   It is bounded as a function of $Z$ since $\PP^n$ is compact.   Looking at the right-hand side again, this implies that $w$ is continuous and weakly admissible.  
(Note that if $w$ is {admissible}, then $\langle\cdot,\cdot\rangle_W$ vanishes on all lines over $H_{\infty}$.)

We also have a notion of sup norm on sections of a line bundle.  
\begin{definition}\rm
Let $\langle\cdot,\cdot\rangle_{W}$ be a Hermitian metric on a line bundle $L$ over $\PP^n$.  Then for each $s\in H^0(\PP^n,L)$ we define
$$
\|s\|_{W}^2 := \sup_{Z\in\PP^n} \left|{\langle s(Z),s(Z)\rangle_{W}} \right|.
$$
\end{definition}

\begin{remark}\rm
Formally, the sup norm $\|p\|_{K,w,k}$ of Definition \ref{def:27} on $\CC[z]_{\leq k}$ can be put into this geometric context  by defining 
$$
W(Z) := \left\{ \begin{matrix}  w(Z/Z_0)^k & \hbox{if } Z_0\neq 0 \hbox{ and } z\in K \\ 
+\infty & \hbox{if } Z_0=0 \hbox{ or } z\not\in K
\end{matrix}\right. 
$$
and defining $\langle\cdot,\cdot\rangle_W$ on $\calO(k)$ as in equation (\ref{eqn:23}).  Then $\langle p(Z),p(Z)\rangle_W=\|p\|_{K,w,k}$.  Note that since $W$ is not necessarily continuous, it is an instance of a more general object called a  \emph{singular Hermitian metric}.  Such objects are important in the application of pluripotential theory to complex geometry \cite{demailly:singular}.
\end{remark}

All of the above goes through on a smooth subvariety $V\subset\CC^n$.  One can define a weight function on $K\subset V$, as well as a Hermitian metric on $\calO(k)$ over $V_{\PP}$.  One simply restricts attention to points of $V$.  

Weight functions are also convenient for doing local computations on a line bundle.  In this paper, we are really only interested in the special case of projective space.

 \begin{example} \label{ex:34}  \rm
  Let $z=(z_1,\ldots,z_n)$ denote affine coordinate in $\CC^n$ and suppose $K\subset V\setminus\{z_1=0\}$, where $V\subset\CC^n$ is an algebraic subvariety, extended to $V_{\PP}\subset\PP^n$.  Let $v=(v_0,v_2,\ldots,v_n)$ be the local coordinates at infinity given by dehomogenization at $z_1$.   This is the holomorphic map $v=g_{10}(z)$ on $\CC^n\cap\{z_1\neq 0\}$  given explicitly by 
$$v_0=\frac{1}{z_1},\ v_2=\frac{z_2}{z_1}, \ \ldots \ , v_n=\frac{z_n}{z_1},$$
and $z=g_{01}(v) := g_{10}^{-1}(v)$ is given by a similar formula.   A section in $H^0(V_{\PP},\calO(k))$ is given by a homogeneous polynomial $p(Z_0,\ldots,Z_n)$, with local evaluations related by 
 $$
 p(1,z_1,\ldots,z_n) = v_0^kp(v_0,1,v_2,\ldots,v_n).
 $$
Hence polynomial evaluation in affine coordinates with weight $w(z)=w(z_1,\ldots,z_n)$ on $K$ transforms to a polynomial evaluation with weight $\tilde w(v) = {v_0}^kw(g_{01}(v))$.  
The transition function simply appears as an additional factor.  This is why it is convenient to allow complex-valued weights.
 \end{example}

\begin{remark}\rm 
Given a positive finite measure $\mu$ supported on $K$ and $k\in\NN$, we also have the \emph{weighted $L^2$ inner product and norm} on $\CC[z]_{\leq k}$, 
$$
\langle p,q\rangle_{\mu,w,k} := \int p(z)_{w,k}\overline{q(z)_{w,k}} d\mu(z), \quad \|p\|_{\mu,w,k}^2 = \langle p,p\rangle_{\mu,w,k} .
$$
The triple $(K,\mu,w)$ is said to satisfy the \emph{Bernstein-Markov property} if there is a sequence $M_1,M_2,\ldots$ of positive integers such that 
$$\|p\|_{K,w,k} \leq M_k\|p\|_{\mu,w,k} \hbox{ for all } k\in\NN,\ p\in\CC[z]_{\leq k}, \quad \hbox{and }   \limsup_{k\to\infty} (M_k)^{1/k} = 1.$$
We call $\mu$ a Bernstein-Markov measure for $K$ (with weight $w$).
\end{remark}

The Bernstein-Markov property is important because it means that certain asymptotic quantities in pluripotential theory associated to a set may be computed using the $L^2$ norm of a Bernstein-Markov measure rather than the sup norm.  The additional tools provided by the $L^2$ theory are important in pluripotential theory, but we will not need them in this paper.





\section{Okounkov bodies and computational algebraic geometry} \label{sec:ok}

We want to study Okounkov bodies associated to varieties using methods of computational algebraic geometry. We will work on an algebraic variety $V\subset\CC^n$. 
 We first review some background material on Noether normalization and normal forms of polynomials. 

\subsection{Normal forms and Noether normalization} 
By the Nullstellensatz, restricting the evaluation of $p\in\CC[z]$  to points of $V$ is equivalent to taking the quotient $\CC[z]/\bI(V)$, with associated equivalence relation $p\sim q$ if $p-q\in\bI(V)$. 
\begin{notation}\rm
Given $k\in\NN$, denote by $\CC[V]_{\leq k}$ the quotient space $\CC[z]_{\leq k}/\sim$ with $\sim$ as above.  For $q\in\CC[z]_{\leq k}$ we can identify equivalence classes containing $q$ under the natural inclusion 
$$
\CC[V]_{\leq k}\ni\{p\in\CC[z]_{\leq k}: p\simeq q\in\CC[V]_{\leq k}\}\hookrightarrow \{p\in\CC[z]: p\simeq q\in\CC[V]_{\leq k}\} \in\CC[V].$$  
Then under this identification, one can see that $\CC[V]=\bigcup_k\CC[V]_{\leq k}$.  
For a general polynomial $p$, put $\deg_V(p)=k$ if $p$ is equivalent to a polynomial of degree $k$ but not of degree $k-1$ (i.e., in $\CC[V]_{\leq k}\setminus\CC[V]_{\leq k-1}$).
\end{notation}
Via dehomogenization in affine coordinates, $H^0(V_{\PP},\calO(k))$ may be identified with $\CC[V]_{\leq k}$.  
  
  \begin{theorem}[Noether Normalization Theorem] \label{thm:noether}
Suppose $V$ is of dimension $m$.  There is a complex linear change of coordinates on $\CC^n$ such that, in the new coordinates (which we denote by  $(x,y):=(x_1,\ldots,x_m,y_1,\ldots,y_{n-m})$), 
\begin{enumerate}
\item The projection map $\pi:V\to\CC^m$ given by $\pi(x,y)=x$ is onto, and $\pi^{-1}(x)$ is finite for each $x\in\CC^m$; 
\item We have an injection $\CC[x]\hookrightarrow\CC[V]$ that exhibits $\CC[V]$ as a finite dimensional algebra over $\CC[x]$.
\end{enumerate}
\end{theorem}
The map $\CC[x]\hookrightarrow\CC[V]$ given in the theorem is given by identifying $p$ with its equivalence class in $\CC[V]$.  (See e.g. Chapter 5 \S 6 of \cite{coxlittleoshea:ideals} for a proof of this theorem.)
  
We turn to algebraic computation in $\CC[V]$; this requires an ordering on monomials.  First, we recall the \emph{lexicographic (lex) ordering} on $\ZZ_{\geq 0}^n$ (denoted $\prec_l$).  We have 
$\alpha\prec_l\beta$ if there exists a $j\in\{1,\ldots,n\}$ for which $\alpha_j<\beta_j$, and $\alpha_k=0$ for all $k>j$.  Monomials in $\CC[z]=\CC[z_1,\ldots,z_n]$ are ordered accordingly: $z^{\alpha}\prec_lz^{\beta}$ if $\alpha\prec_l\beta$, so that $z_1\prec_l z_2\prec_l\cdots\prec_l z_n$.  We will come back to lex ordering later.  

We also recall the \emph{grevlex ordering} which has good computational properties.  This is the ordering $\prec_g$ for which $z^{\alpha}\prec_g z^{\alpha'}$ whenever
\begin{enumerate}
\item $|\alpha|<|\alpha'|$; or
\item $|\alpha|=|\alpha'|$ and $z^{\alpha}\prec_l z^{\alpha'}$.
\end{enumerate}

\begin{notation}\rm 
For a polynomial $p\in\CC[z]$l, let us denote by $\lt(p)$ the leading term of $p$ with respect to grevlex, and for an ideal $I$ of $\CC[z]$, let 
$
\lt(I):=\{\lt(p): \ p\in I\}.$
 \end{notation}
 
We will use the grevlex ordering to compute {normal forms}. Let us recall what these are.  First, a \emph{Groebner basis} of an ideal $I$ is a collection $\{g_1,\ldots,g_{\ell}\}\subset I$ for which
 $$
 I=\langle g_1,\ldots,g_{\ell}\rangle \hbox{ and } \langle\lt(I)\rangle = \langle \lt(g_1),\ldots,\lt(g_{\ell})\rangle.
 $$ 
 For each element of $\CC[z]/I$ there is a unique polynomial representative, called the \emph{normal form}, which contains no monomials in the ideal $\langle\lt(I)\rangle$.  The normal form of a polynomial $p$ may be computed in practice as the remainder $r$ upon dividing $p$ by a Groebner basis of $I$:
 $$
 p = q_1g_1+\cdots+q_{\ell}g_{\ell} + r,
 $$
 where $q_1,\ldots,q_{\ell}\in\CC[z]$ are the quotients.  (See e.g. chapter 3 of \cite{coxlittleoshea:ideals} for a description of the associated division algorithm.)
 
 Let $\CC[z]_I$ be the collection of normal forms.  This is an algebra over $\CC[z]$ under the usual addition of polynomials, and with multiplication defined by 
 \begin{equation}\label{eqn:nf} \CC[z]_I\times\CC[z]_I\ni (r_1,r_2)\longmapsto \hbox{``the normal form of $r_1r_2$''}\in\CC[z]_I.\end{equation}

The following algebraic version of Noether normalization is given in \cite{coxmau:transfinite}.  

\begin{proposition} \label{prop:36}
Let $V$ be of dimension $m$ and let $(x,y)$ be coordinates as in Theorem \ref{thm:noether}.  Let $\CC[x,y]_{I}$ be the algebra of normal forms for $I=\bI(V)$.   Then 
\begin{enumerate}
\item \label{prop3.6(1)} $\deg_V(p)=\deg(p)$ (i.e. the usual degree) whenever $p$ is a normal form.
\item \label{prop3.6(2)} We have the inclusion $\CC[x]\subseteq\CC[x,y]_I$, which exhibits $\CC[x,y]_I$ as a finite dimensional algebra over $\CC[x]$.  
\end{enumerate}
\end{proposition}
 Here, multiplication in $\CC[x,y]_I$ is as in equation (\ref{eqn:nf}).  The proposition says that any $p\in\CC[x]$ is a normal form, and shows that grevlex has good computational properties.  
 
  In what follows, we will usually assume polynomials to be normal forms, and $\CC[V]$ will be identified with  $\CC[x,y]_{I}$.   The inclusion in item (\ref{prop3.6(2)}) of the proposition is called a \emph{Noether normalization}; we will also write (via our identifications) $\CC[x]\subseteq\CC[V]$.  Let us also refer to the coordinates $(x,y)$ as \emph{(Noether) normalized coordinates}.
  
  Since $\CC[V]$ is finite dimensional over $\CC[x]$, and has a basis of monomials, there are only a finite number of monomials $y^{\beta}$ for which $x^{\alpha}y^{\beta}$ is a normal form.  Hence any normal form, being a linear combination of such monomials, can be expressed as a finite sum
 \begin{equation} \label{eqn:normalform}
 p(z) = p(x,y) \  = \  \sum_{\beta} y^{\beta}p_{\beta}(x), \quad p_{\beta}\in\CC[x].  
 \end{equation}
 
 \begin{example} \label{ex:37} \rm
 The (complexified) sphere in $\CC^3$ is given by
 \begin{equation}\label{eqn:0.1} V=\{z=(z_1,z_2,z_3)\in\CC^3: z_1^2+z_2^2+z_3^2=1\},\end{equation}
and $\langle z_1^2+z_2^2+z_3^2-1\rangle=\bI(V)=:I$.  Any polynomial in $I$ is of the form $q(z)(z_1^2+z_2^2+z_3^2-1)$, and
 $$\lt(q(z)(z_1^2+z_2^2+z_3^2-1))=\lt(q(z)z_3^2)\in\langle z_3^2\rangle;$$ it follows easily that $\langle\lt(I)\rangle=\langle z_3^2\rangle$.  Hence a normal form $p\in\CC[V]$ is a polynomial given by 
 \begin{equation} \label{eqn:normalformsphere}
 p(z)= p_1(z_1,z_2) + z_3p_2(z_1,z_2), \quad p_1,p_2\in\CC[z_1,z_2].
 \end{equation}
(Compare the above to (\ref{eqn:normalform}).)  As a $2$-dimensional algebra over $\CC[z_1,z_2]$, multiplication is given by 
 $$
 (p_1+z_3p_2)\cdot(q_1+z_3q_2) \ = \ p_1q_1+(1-z_1^2-z_2^2)p_2q_2 \  + \  z_3(p_1q_2+p_2q_1).
 $$
Clearly, $x=(x_1,x_2):=(z_1,z_2)$ and $y=z_3$ give Noether normalized coordinates satisfying Theorem \ref{thm:noether}: $x\in\CC^2$ lifts to at most 2 points $(x,y)\in V\subset\CC^3$ given by the branches of the square root in the expression  $y=(1-x_1^2-x_2^2)^{1/2}$.  
 \end{example}

 \subsection{Okounkov body computation}
 
 Following Nystr\"{o}m \cite{nystrom:transforming}, let us define the Okounkov body of a line bundle.  Returning to the geometric setting, let $L$ be a holomorphic line bundle over a complex manifold $M$ of dimension $n$, and $p\in M$.  In a local trivialization containing $p$, any $s\in H^0(M,L)$ is given by a holomorphic function (let us also denote the function by $s$).  Hence it can be expressed as a power series
\begin{equation}\label{eqn:31}
s = \sum_{\alpha} c_{\alpha}z^{\alpha}
\end{equation}
where $z$ is a local holomorphic coordinate centered at the point $p$, and we use multi-index notation: for $\alpha=(\alpha_1,\ldots,\alpha_n)\in\ZZ_{\geq 0}^n$, we have  
$z^{\alpha}=z_1^{\alpha_1}\cdots z_n^{\alpha_n}$.

\begin{definition}\rm
Suppose a local holomorphic coordinate $z$ is fixed at $p\in M$.  Given a section $s\in H^0(M,L)$, we define $\nu(s)$ to be the lowest exponent in the power series (\ref{eqn:31}) with respect to the lex order, $\prec_l$.  That is, if $\nu(s)=\gamma$ then $c_{\alpha}=0$ whenever $\alpha\prec_l\gamma$.  When $\gamma=\nu(s)$ we also define
\begin{itemize}
\item[]$\ttt(s)=c_{\gamma}z^{\gamma}$, the \emph{trailing term};
\item[]$\tc(s)=c_{\gamma}$, the \emph{trailing coefficient}; and
\item[]$\tm(s)=z^{\gamma}$, the \emph{trailing monomial}.
\end{itemize}
\end{definition}
 
\begin{definition}\rm
Fix a local holomorphic coordinate $z$ at a point $a\in M$, and let $k\in\NN$. Expand $s\in H^0(M,L^{\otimes k})$ as in (\ref{eqn:31}), and define   
$$
\calN_k:= \{\nu(s)\in\ZZ_{\geq 0}^n: s\in H^0(M,L^{\otimes k})   \}.
$$
Let $\Delta_k\subset\RR^m$ be the convex hull of the set $\frac{1}{k}\calN_k$.  The \emph{Okounkov body of $L$ (with respect to these coordinates)}, denoted by $\Delta=\Delta(L)$,  is defined to be the convex hull of the set $\bigcup_{k\in\NN} \Delta_k$.
\end{definition}

In our concrete setting of an algebraic subvariety $V\subseteq\CC^n$ of dimension $m$, we will use the sections $H^0(V_{\PP},\calO(k))$, or equivalently, the polynomials $\CC[V]_{\leq k}$.  (For convenience, let us assume $V_{\PP}$ is smooth, so that the theory on a complex manifolds can be transferred without any technicality.) 
  We will use Noether normalized coordinates $(x,y)=(x_1,\ldots,x_m,y_1,\ldots,y_{n-m})$ on $V$.  Without loss of generality, assume that the point $a$ at which the Okounkov body is calculated is of the form $(0,\ldots,0,a_{m+1},\ldots,a_n)$, i.e., all of the $x$ coordinates are zero.\footnote{Any translation $x\mapsto x+c=\tilde x$ gives an isomorphism of normal forms $p(\tilde x,y)\mapsto p(x+c,y)$; the verification of this is left as an exercise.}  We will also assume that the point $a$ is a regular point for the projection to the $x$ coordinates, i.e., $V$ satisfies the hypotheses of the holomorphic implicit function theorem:
 \begin{equation}\label{eqn:implicit}
\det\begin{bmatrix} \partial f_1/\partial y_1(a) & \cdots & \partial f_1/\partial y_{n-m}(a) \\
\vdots & \ddots & \vdots \\
\partial f_{n-m}/\partial y_1(a) & \cdots & \partial f_{n-m}/\partial y_{n-m}(a)
\end{bmatrix} \neq 0,
 \end{equation}
 where the polynomials $f_1,\ldots, f_{n-m}$ determine $V$ in a neighborhood of $a$.  Locally, we can write $y=Y(x)$ for some holomorphic function $Y$ in a neighborhood of $x=(0,\ldots,0)$.  The $x$ coordinates in the Noether normalization provide the local coordinates with which the Okounkov body will be calculated.  
  \begin{notation}\rm
 Suppose $Y(x)=(Y_1(x),\ldots,Y_{n-m}(x))$ in components.  Then for a multi-index $\beta=(\beta_1,\ldots,\beta_{n-m})$ the holomorphic function $Y^{\beta}$ is given by
 $$Y^{\beta}(x) := Y_1(x)^{\beta_1}Y_2(x)^{\beta_2}\cdots  Y_{n-m}(x)^{\beta_{n-m}}.
 $$  In a neighborhood of the origin,  it may be expressed as a power series in $x$:
 \begin{equation}\label{eqn:Yps}
 Y^{\beta}(x)= \sum_{\alpha} c_{\beta\alpha}x^{\alpha}.
 \end{equation}
 \end{notation}
  
Let $k\in\NN$; we want to compute $\calN_k\subset\ZZ_{\geq 0}^m$.  A section of $H^0(V_{\PP},\calO(k))$ can be identified with a polynomial in $\CC[V]_{\leq k}$; by Proposition \ref{prop:36}(\ref{prop3.6(1)}) it is given by a normal form $p$ of degree $\leq k$.  Let $\{y^{\beta}\}_{\beta}$ be the finite collection of $y$ monomials as in (\ref{eqn:normalform}).  In a neighborhood of the origin, we rewrite $p=p(x,y)$ as the holomorphic function $x\mapsto p(x,Y(x))$.  Write this as a power series in $x$, and compute the coefficients by multiplying out the terms of the polynomials $p_{\beta}$ with the power series (\ref{eqn:Yps}) for $Y^{\beta}$:
$$
p(x,Y(x)) = \sum_{\beta} Y^{\beta}(x)p_{\beta}(x) = \sum_{\beta}\sum_{\alpha} p_{\beta}(x)c_{\beta\alpha}x^{\alpha} =: \sum_{\alpha} b_{\alpha}x^{\alpha}.
$$
We can then read off $\nu(p)$ from the trailing term on the right-hand side.

\medskip

 We determine the possibilities for $\nu(p)$.  For all $\alpha\in\ZZ_{\geq 0}^m$ with $|\alpha|\leq k$, we have, by Proposition \ref{prop:36},
$$
x^{\alpha}\in\CC[x]_{\leq k}\subseteq\CC[V]_{\leq k}, \quad\hbox{and } \nu(x^{\alpha}) = \alpha.
$$
The points $\alpha/k$, with $\alpha$ as above, fill out a grid of rational points in the region bounded by the coordinate hyperplanes and the standard simplex in $\RR^m$.\footnote{The standard simplex in $\RR^m$ is the convex hull of the standard basis vectors $\{(1,0,\ldots,0),\ldots,(0,\ldots,0,1)\}$.   }  This accounts for all normal forms in $\CC[x]_{\leq k}$, and so remaining points in $\calN_k$ must be calculated from normal forms containing powers of $y$, which involve the analytic functions $Y^{\beta}$.  We will use the notion of \emph{$S$-polynomial} to do this systematically (Definition \ref{def:Spoly}).  Computing the points of $\frac{1}{k}\calN_k$ for larger and larger values of $k$, we eventually fill in the Okounkov body.

We write out the details of some calculations on the complexified sphere in what follows.  From these calculations, we derive a general method that can be applied to varieties given in Noether normalized coordinates.\footnote{But we do not give a rigorous proof.}

\subsection{Explicit computation on the sphere}
Let $V$ be the complexified sphere as in (\ref{eqn:0.1}), Example \ref{ex:37}.  For the Okounkov body, take the point $a=(0,0,1)$, and $(z_1,z_2)$ as the local coordinates in which to expand polynomials on $V$.  By (\ref{eqn:normalformsphere}) the normal forms in $\CC[V]$ are linear combinations of monomials of the form $z_1^jz_2^k$ or $z_1^jz_2^kz_3$.  At the point $a$, 
$$\frac{\partial}{\partial z_3}(z_1^2+z_2^2+z_3^2-1)\bigl|_{(0,0,1)}=2z_3\bigl|_{(0,0,1)}=2\neq 0,$$ 
so (\ref{eqn:implicit}) holds and we can write $z_3=Y(z_1,z_2)$, which is in fact the standard square root function:
$$
z_3 = Y(z_1,z_2) = (1-z_1^2-z_2^2)^{1/2} =: \sum c_{jk}z_1^jz_2^k.
$$
The coefficients $\displaystyle c_{jk}=\frac{\partial^{j+k}}{\partial z_1^j\partial z_2^k}Y(0,0)$ may be computed by implicit differentiation, for example,
\begin{equation*}
0 = \frac{\partial}{\partial z_1}(z_1^2+z_2^2+Y^2)\Bigl|_{(0,0,1)}  
= 2z_1 + 2Y\frac{\partial Y}{\partial z_1}\Bigl|_{(0,0,1)} = 2\cdot 0 + 2c_{10}, 
\end{equation*}
so that $c_{10}=0$.  Further differentiation gives more coefficients, e.g. $c_{20}=-\frac{1}{2}$, so that $z_3=1-\frac{1}{2}z_1^2+\cdots$.  (Note that here it is faster just to read off the coefficients from the binomial series for the square root.)

The following properties are straightforward to verify.
\begin{lemma}\label{lem:39b}
For all nonzero $p,q\in\CC[V]$, 
\begin{eqnarray}
\nu(pq) &=& \nu(p)+\nu(q), \hbox{ and}  \nonumber   \\
\nu(p+q) &\succ_l&  \min\{\nu(p),\nu(q)\}.   \label{eqn:nu}
\end{eqnarray}
When $\nu(p)\neq \nu(q)$ then $\nu(p+q)=\min\{\nu(p),\nu(q)\}$. \qed
\end{lemma} 
Hence (\ref{eqn:nu}) is strict only if $\nu(p)=\nu(q)$ and we have cancellations of lowest terms.  This motivates the following definition.
\begin{definition}  \label{def:Spoly}\rm
Suppose $\nu(p)=\nu(q)$.  We define the \emph{$S$-polynomial\footnote{Here $S$ stands for \emph{syzygy} (a pair of connected or corresponding things).  The terminology is adapted from chapter 2 of \cite{coxlittleoshea:ideals}.}  of $p$ and $q$} by  
$$S(p,q):= \tc(q)p - \tc(p)q.$$
\end{definition}

From $p=a_0+a_1z_1+a_2z_2+a_3z_3\in\CC[V]_{\leq 1}$, the properties of $\nu$ under addition in Lemma \ref{lem:39b} gives possible values for $\nu(p)\in\calN_1$ as
$$
\nu(1),\nu(z_1),\nu(z_2),\nu(z_3).
$$
Now $\nu(1)=(0,0)=\nu(z_3)$, so we can arrange a possible cancellation of the constant term. We compute  $S(z_3,1) = z_3-1$ and $\nu(z_3-1)=(2,0)$, to get the remaining point of $\calN_1$. 

 The following elementary observation is useful.
\begin{lemma} \label{lem:212}
For any $k\in\NN$, $dim\CC[V]_{\leq k}$ is the number of points in $\calN_k$.
\end{lemma}
\begin{proof}
Let $M_k:=\dim\CC[V]_{\leq k}$ and let $\{\be_j\}_{j=1}^{M_k}$ be a basis.  Then $\{\nu(\be_j)\}_{j=1}^{M_k}\subseteq\calN_k$, so the number of points in $\calN_k$ is at least $M_k$.  

On the other hand, let $p,q$ be nonzero polynomials in $\CC[V]_{\leq k}$ with $\nu(p)\neq\nu(q)$; without loss of generality, $\nu(p)\prec_l\nu(q)$.  Then 
$$\nu(\lambda p+\mu q)= \left\{\begin{matrix}  \nu(p)  & \hbox{if }\lambda\neq 0 \\
\nu(q) & \hbox{if } \lambda=0\neq\mu 
\end{matrix}  \right.  .$$
This implies $\lambda p+\mu q\neq 0$ unless $\lambda=\mu=0$.  Hence any set of polynomials $\{p_{\alpha}\}_{\alpha\in\calN_k}\subset\CC[V]_{\leq k}$ (for which $\nu(p_{\alpha})=\alpha$) is linearly independent, so the size of this set is at most $\dim\CC[V]_{\leq k}$.
\end{proof}

Let us compute $\calN_2$; this will motivate the general case.  First, note that products of pairs of polynomials from $\calB_1:=\{1,z_1,z_2,z_3,S(z_3,1)\}$ span $\CC[V]_{\leq 2}$; this follows easily from the fact that the map
\begin{equation} \label{eqn:37a}
\CC[v_1,v_2,v_3,v_4]_{\leq 2} \  {\longrightarrow} \  \CC[V]_{\leq 2}
\end{equation}
given by making the substitutions $v_1=z_1,v_2=z_2,v_3=z_3,v_4=S(z_3,1)$,  followed by reduction to normal form,  is well-defined and onto.  From such products, we immediately obtain points of $\calN_2$ given by
$$\nu(pq), \quad \hbox{where }p,q\in\calB_1.$$
This gives nine points.  Since $\dim\CC[V]_{\leq 2}=9$, then by the previous lemma, we are done.

Without knowing the dimension a priori, one can simply compute possible $S$-polynomials and check that no new points are obtained.  For $\calN_2$,  we have $\nu(z_1^2)=\nu(z_3-1)=(2,0)$, and 
\begin{eqnarray*}
S(S(z_3-1),z_1^2))=z_3-1+\tfrac{1}{2}z_1^2
\end{eqnarray*}
with $\nu(z_3-1+\tfrac{1}{2}z_1^2)=(4,0)$.  But now, $\nu(S(z_3,1)^2)=(4,0)$ also, and provides another possibility for cancellation.  We calculate $(z_3-1)^2=z_3^2-2z_3+1 = 2z_3-2-z_1^2-z_2^2$, and then the $S$-polynomial 
$$S(z_3-1+\tfrac{1}{2}z_1^2,2z_3-2-z_1^2-z_2^2) = \tfrac{1}{8}z_2^2,$$ whose $\nu$ value is $(2,0)$. Finally, $S(\frac{1}{8}z_2^2,z_2^2)=0$.  There are no further $S$-polynomials to be calculated.  See Figure 1 for a picture.

\begin{figure} \label{fig0}
\includegraphics[height=100mm]{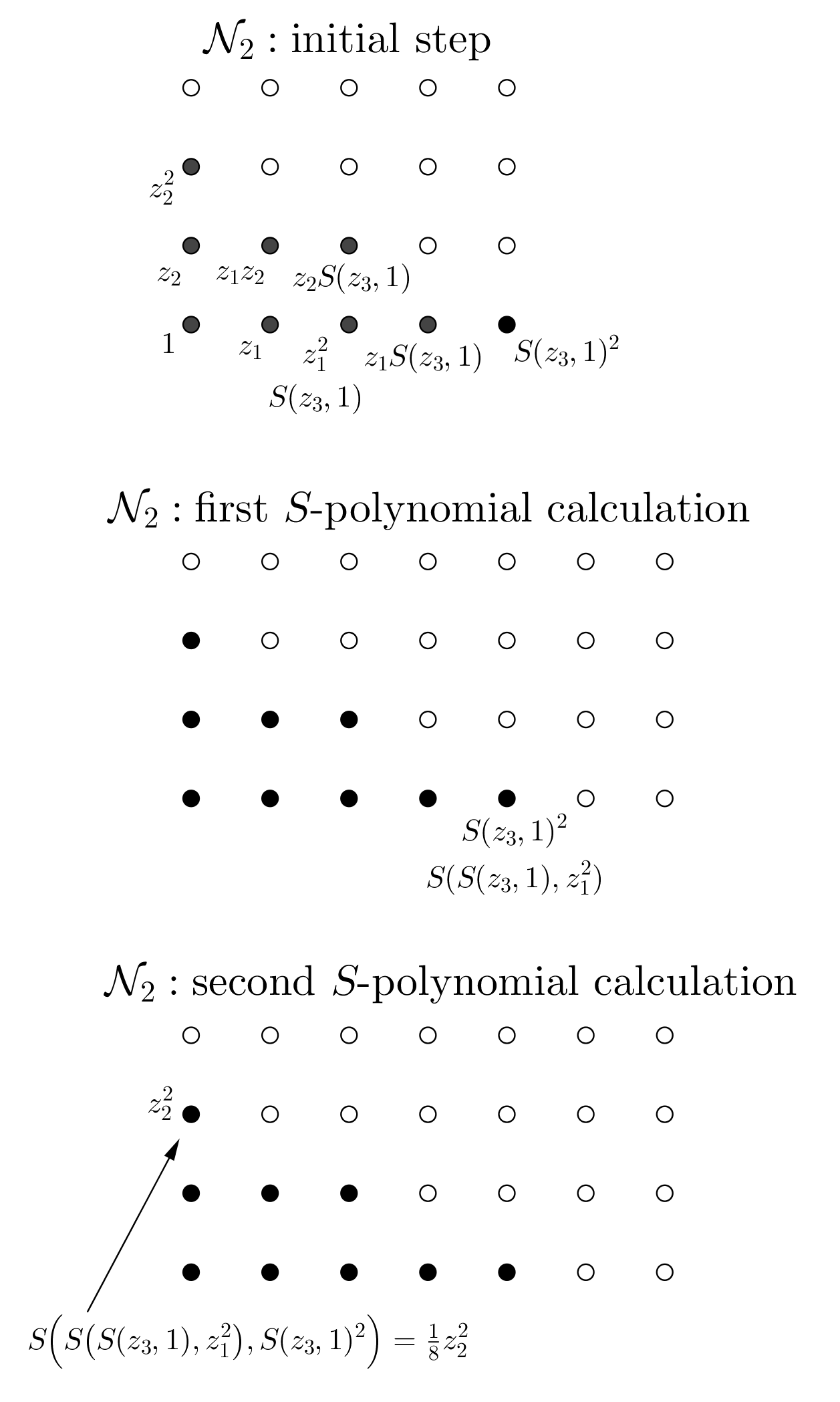}
\caption{Computing points of $\calN_2$ and $S$-polynomials.}
\end{figure}

\subsection{General algorithm} Let us use the same method inductively for $\calN_k$.  We begin by setting $\calN_k$ to be $\calN_{k-1}$.  Step 1 below uses the fact that if the set $\{\be_j\}$ spans $\CC[V]_{\leq k-1}$, then the normal forms of the products $$\{\be_j\be_{\ell}: \deg(\be_j\be_{\ell})\leq k  \}$$ span $\CC[V]_{\leq k}$.   
\begin{itemize}
\item[Step 1] {\it Multiplication}.  
Let $\{\be_j\}$ be a collection of polynomials whose $\nu$ values give $\calN_{k-1}$ (one polynomial for each $\nu$ value).  
Compute all products in which $\deg(\be_j\be_{\ell})=k$, and denote this collection of products by $\calP$.

\item[Step 2] {\it Compute $\nu$ values}. Compute $\nu(p)$ for each $p\in\calP$.  Add any new $\nu$ values thus obtained to $\calN_k$.  

\item[Step 3] {\it Find $\alpha$ for possible cancellation}.  Let $\alpha\in\calN_k$ be the smallest (according to $\prec_l$) point that is given by more than one element of $\{\be_j\}\cup\calP$.

\item[Step 4]  {\it Compute $S$-polynomials at $\alpha$}. Compute $S$-polynomials associated to distinct pairs of polynomials $p$ in $\{\be_j\}\cup\calP$ for which $\nu(p)=\alpha$.

\item[Step 5] {\it Update $\calN_k$ and $\calP$}.  Find the $\nu$ values of the nonzero $S$-polynomials in the previous step.  Add the $\nu$ values to the set $\calN_k$, and add the $S$-polynomials to the set  $\calP$.

\item[Step 6] {\it Repeat}.  Let $\alpha$ be the next point in $\calN_k$ (according to $\prec_l$) that is the $\nu$ value of at least two different elements of $\calP\cup\{\be_j\}$.  But if no such $\alpha$ exists, or we know that the current size of $\calN_k$ equals the dimension of $\CC[V]_{\leq k}$, go to the next step. Otherwise, return to Step 4.

\item[Step 7]  {\it Stop}.  Output $\calN_k$.
\end{itemize}

To compute $\calN_{k+1}$, expand the collection $\{\be_j\}$ by adding elements of $\calP$ computed in the above algorithm so that there is one polynomial for each point of $\calN_k$; then begin again from Step 1 with $k:=k+1$.

The above method can be applied to a subvariety $V\subseteq\CC^n$ of pure dimension $m$ at a smooth point $a\in V$.    Assume that we are working in coordinates $(x_1,\ldots,x_m,y_1,\ldots,y_{n-m})=(x,y)$ that give a Noether normalization of $V$ at $a$. Without loss of generality, we translate  coordinates so that the $x$ coordinates of $a$ are zero.  

  Start with the monomials $\{1,x_1,\ldots,x_m,y_1,\ldots,y_{n-m}\}$ and compute $\calN_1$ using the power series of the monomials $y_j$, $j=1,\ldots,n-m$, as local functions of the $x$ coordinates near $a$.   Continue inductively, computing $\calN_2,\calN_3,\ldots$ exactly as above, by multiplying polynomials, and computing $S$-polynomials if required.  Use the points of $\frac{1}{k}\calN_k$ for sufficiently many values of $k$ to fill in the Okounkov body.  
  
  \begin{remark}\rm
The above algorithm is analogous to Gaussian elimination. Eliminating a lower order monomial by computing $S$-polynomials is essentially  eliminating a variable in a linear system.
\end{remark}

Using Lemma \ref{lem:212} it may be possible to deduce the Okounkov body exactly.  Let us do this for our sphere example.
\begin{proposition} \label{prop:oksphere}
The Okounkov body of the complexified sphere is the triangle given by 
$$
\Delta :=\{(\theta_1,\theta_2)\in\RR^2: \theta_1,\theta_2\geq 0, \theta_1+2\theta_2\leq 1\}.
$$
\end{proposition}

\begin{proof}
It is sufficient to show that for each $k\in\NN$, $\Delta\subseteq\Delta_k$.  This will follow from the following claim: \\ 
\centerline{\sl $\calN_k$ consists of the integer points of $k\Delta$, i.e., $\calN_k=k\Delta\cap\ZZ^2$.} \\
  Our calculations above have already established this claim for $k=1,2$.  Working by induction on $k$, suppose $\calN_{k-1}=(k-1)\Delta\cap\ZZ^2$.    Then it is easy to see that 
$$k\Delta\cap\ZZ^2= \{\nu(pq):\ p\in\CC[V]_{\leq 1},\ q\in\CC[V]_{\leq k-1}   \} \subseteq\calN_k.$$ 
The claim will follow if the number of points in $k\Delta\cap\ZZ^2=:\calM_k$ equals $\dim\CC[V]_{\leq k}$.  By induction, this can be reduced to verifying that the number of new points, i.e., the size of $\calM_k\setminus\calN_{k-1}$, is equal to the number of elements of degree $k$ in a basis for $\CC[V]_{\leq k}$. 

On the one hand, the elements of degree $k$ in a monomial basis of $\CC[V]_{\leq k}$ are $$\{z_1^k,z_1^{k-1}z_2,\ldots,z_2^k,z_1^{k-1}z_3,z_1^{k-2}z_2z_3,\ldots,z_2^{k-1}z_3\},$$ which are $2k-1$ elements.  On the other hand, one can see that $\calM_k\setminus\calN_{k-1}$ has $2k-1$ points by some elementary geometry in the plane.\footnote{Precisely, it consists of the $k$ points along the hypotenuse and the $k-1$ points on the neighbouring parallel line.}
\end{proof}

We close this section by looking a bit more closely at the limiting behaviour of the Okounkov body construction.  Although $k<l$ does not imply $\Delta_k\subset\Delta_l$, it is easy to verify the weaker statement
\begin{equation}\label{eqn:mon}
\Delta_k\subseteq \Delta_{rk}, \quad\hbox{for all } k,r\in\NN.
\end{equation}
This follows from the fact that if $\alpha\in\frac{1}{k}\calN_k$ then $k\alpha=\nu(p)$ for some $p\in\CC[V]_{\leq k}$.  Consequently, $p^r\in\CC[V]_{\leq rk}$ with $\nu(p^r)=rk\alpha$, so $\alpha\in\frac{1}{rk}\calN_{rk}$.  Hence $\frac{1}{k}\calN_k\subset\frac{1}{rk}\calN_{rk}$,  and (\ref{eqn:mon}) follows upon taking convex hulls. 

The following result will be used in the next section.  

\begin{proposition} \label{prop:215a}
Let $\Delta\subset\RR^m$ be the Okounkov body of some variety, $\Delta^{\circ}$ its interior, and  $Q\subset\Delta^{\circ}$ a compact convex subset with nonempty interior.  For a sufficiently large positive integer $k_0$,  
\begin{equation}\label{eqn:l215}
\tfrac{1}{k_0}\calN_{k_0}\cap Q = \tfrac{1}{k_0}\ZZ^m\cap Q. 
\end{equation}
If $|\calN_k\cap kQ|$ denotes the number of points in $\calN_k\cap kQ$, then
\begin{equation} \label{eqn:l215b}
\tfrac{1}{k^m}|\calN_k\cap kQ| \longrightarrow \vol(Q) \quad\hbox{as }k\to\infty. 
\end{equation}
\end{proposition}

Although the proof is elementary, it is rather technical, so we omit it.  In Lemma 2.3 of \cite{nystrom:transforming}, it is shown (using a theorem of Khovanskii in \cite{khovanskii:newton}) that there is some constant $C$ such that (\ref{eqn:l215}) holds as soon as the distance from $Q$ to $\partial\Delta$ is greater than $C/k_0$.  

We can give an idea for why (\ref{eqn:l215}) holds.  For $Q$ compact in $\Delta^{\circ}$, the number of points in $kQ\cap\calN_k$ is almost the same as the number of points in $kQ\cap\calN_{k-1}$.  
 By taking certain products and powers of polynomials (as in the verification of (\ref{eqn:mon})), one can fill in the gaps by `translation from existing nodes'. See Figure 2 for an illustration.
 
Equation (\ref{eqn:l215b}) follows easily from (\ref{eqn:l215}).  Up to a boundary correction of order $k^{m-1}$, the number of points in $\tfrac{1}{k_0}\ZZ^m\cap Q$ is equal to the number of $m$-dimensional cubes in the corresponding grid that cover $Q$, and the volume of each cube is $(\frac{1}{k})^m$.  The total volume of these cubes converges to the total volume of $Q$.  

Note that $\Delta$ is compact: for sufficiently large $k$, the number of points in $\calN_k=\dim\CC[V]_{\leq k}$ is given by the Hilbert polynomial of $V$, which is a polynomial in $k$ of degree $m$. Hence $\tfrac{1}{k^m}|\calN_k|$ is uniformly bounded.  The same argument as in the previous paragraph shows that $\tfrac{1}{k^m}|\calN_k|\to\vol(\Delta)$, which must be finite.
 
  \begin{figure} \label{figOk}
\includegraphics[height=50mm]{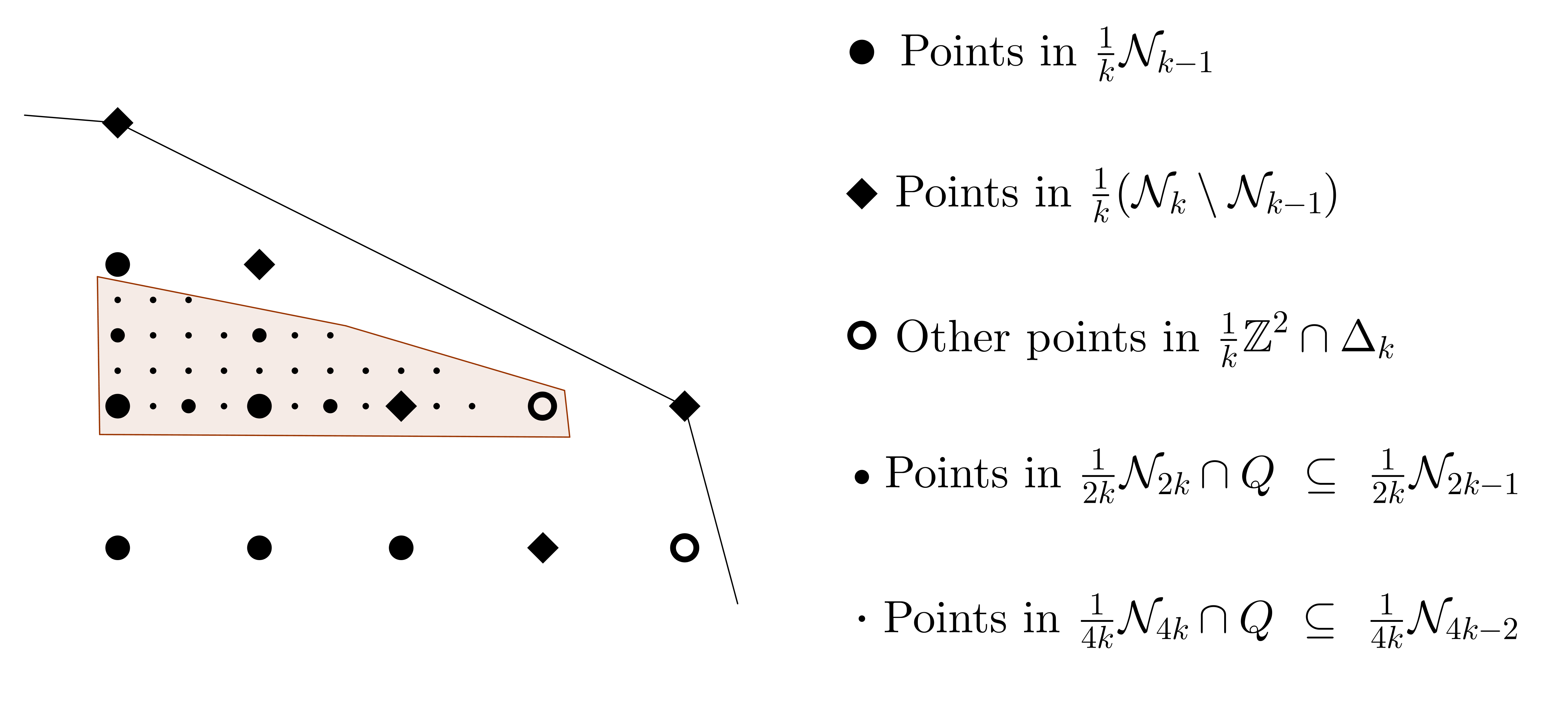}
\caption{ 
Suppose $\frac{1}{k}\nu(p)=\alpha$ where  $p\in\CC[V]_{\leq k-1}$ (one of the large black circles).    
Then $$p^2\in\CC[V]_{\leq 2k-2}  \ \hbox{ and }  \  z_1p^2,z_2p^2\in\CC[V]_{\leq 2k-1},$$ giving $\tfrac{1}{2k}\nu(p^2)=\alpha$, $\tfrac{1}{2k}\nu(z_1p^2)=\alpha+(\tfrac{1}{2k},0)$, and $\tfrac{1}{2k}\nu(z_2p^2)=\alpha+(0,\tfrac{1}{2k})$ as points of $\calN_{2k}$.  We can do the same sort of thing to get points of $\calN_{4k}$.  
In the above example, the `missing point' of $Q$ in $\frac{1}{k}\ZZ^m\setminus\frac{1}{k}\calN_k$ will be covered with a point of $\tfrac{1}{8k}\calN_{8k}$; hence (\ref{eqn:l215}) holds with $k_0=8k$. }
 \end{figure}

\section{Chebyshev transform and transfinite diameter} \label{sec:dir}

Let $V\subset\CC^n$ be an algebraic variety.  Throughout this section we will assume that we are working in local coordinates on $V$ constructed from a Noether normalization, with the associated Okounkov body as defined in the previous section. Let us denote these local coordinates by $z$.   We will define directional Chebyshev constants associated to polynomials on $V$, and the Chebyshev transform.   Then we will define transfinite diameter in terms of a special basis of $\CC[V]$.  Finally, we will prove the main theorem that gives transfinite diameter as an integral of the Chebyshev transform over the Okounkov body.  

  Given $k\in\NN$ and a multi-index $\alpha\in \calN_k$, we define the normalized class of polynomials 
$$
\calM(k,\alpha) := \bigl\{ p\in\CC[V]_{\leq k}: \ttt(p) = z^{\alpha} \hbox{ (i.e., $\tc(p)=1$)} \bigr\}.
$$

Let $K\subset V$ be a set and $w$ an admissible weight on $K$. 

 \addtocounter{footnote}{-4}

\begin{definition}\label{def:Cheby}
\rm For $k$ a positive integer and $\alpha\in \calN_k$.  We define the \emph{Chebyshev constant $T^w_k(K,\alpha)$} by 
$$
T^{w}_k(K,\alpha):= \inf\bigl\{  \|p\|_{K,w,k}: \   p\in\calM(k,\alpha)  \bigr\}^{1/k}.
$$
A polynomial $t$ for which $T^w_k(K,\alpha) = \|t\|_K^{1/k}$ will be called a \emph{Chebyshev polynomial.}
\end{definition}

\begin{lemma} \label{lem:subm}
We have
$$
T^w_{k_1+k_2}(K,\alpha_1+\alpha_2)^{k_1+k_2}   \leq    T^w_{k_1}(K,\alpha_1)^{k_1}T^w_{k_2}(K,\alpha_2)^{k_2}  .
$$
\end{lemma}

\begin{proof}
Let $t_{k_1},t_{k_2}$ satisfy  $\|t_{k_1}\|_{K,w,k_1}=T^w_{k_1}(K,\alpha_1)^{k_1}$ and $\|t_{k_2}\|_{K,w,k_2}=T^w_{k_2}(K,\alpha_2)^{k_2}$.  Then $t_{k_1}t_{k_2}\in\calM(k_1+k_2,\alpha_1+\alpha_2)$, and so
$$
T_{k_1+k_2}^w(K,\alpha_1+\alpha_2) \leq  \|t_{k_1}t_{k_2}\|_{K,w,k_1+k_2} \leq \|t_{k_1}\|_{K,w,k_1}\|t_2\|_{K,w,k_2}.
$$
\end{proof}

The lemma says that $(k,\alpha)\mapsto T^w_k(K,\alpha)^k$ is a \emph{submultiplicative function}.  This property allows us to define directional Chebyshev constants via a limiting process. 

\begin{lemma}\label{lem:clim}
 The limit
$$ 
T^{w}(K,\theta) := \lim_{\substack{k\to\infty \\ \alpha/k\to\theta}} T^w_k(K,\alpha)
$$ 
exists for each $\theta\in \Delta^{\circ}$. 
\end{lemma}

\begin{proof}
Let $\{\alpha_j\},\{\beta_j\},\{s_j\},\{t_j\}$ be sequences of multi-indices such that as $j\to\infty$, we have 
\begin{equation}\label{eqn:sjtj}  s_j,t_j\to\infty, \quad \frac{\alpha_j}{s_j},\frac{\beta_j}{t_j} \to\theta,\end{equation}
  and 
$$
T^w_{s_j}(K,\alpha_j) \to L_1:= \liminf_{\substack{k\to\infty\\ \alpha/k\to\theta}} T^w_k(K,\alpha), \ \hbox{ and } \  
T^w_{t_j}(K,\beta_j) \to L_2:= \limsup_{\substack{k\to\infty\\ \alpha/k\to\theta}} T^w_k(K,\alpha).
$$
The existence of such an approximating sequence follows from Proposition \ref{prop:215a}.\footnote{Apply the proposition with $Q$ being a small closed $m$-dimensional cube containing $\theta$.}  

It is sufficient to show that $L_2\leq L_1$.  To this end, let $\epsilon>0$, and choose an index $j_0$ large enough that $T_{s_{j_0}}(K,\alpha_{j_0})^{s_{j_0}} < L_1-\epsilon$.  Choose a polynomial $p\in\calM(s_{j_0},\alpha_{j_0})$ such that $\|p\|_{K,w,s_{j_0}}^{1/s_{j_0}}<L_1-\epsilon$.  Now, for each $j$ let $q_j$ be the largest non-negative integer for which
$$
\beta_j = q_j\alpha_{j_0} + \ell_j 
$$
and all entries of $\ell_j$ are non-negative.   Using the fact that all entries of $\alpha_{j_0}$ are positive, and each entry of $\beta_j$ goes to infinity, it is straightforward to verify that $q_j\to\infty$ and the set of $\ell_j$'s is a finite subset of $\ZZ^m$.  Using (\ref{eqn:sjtj}), we also have $\frac{q_js_{j_0}}{t_j}\to 1$ as $j\to\infty$.  

For every $j$, consider the polynomial $P_j:=p^{q_j}r_j$, where $r_j$ is a polynomial for which $\nu(r_j)=\ell_j$.  Then $P_j\in\calM(t_j,\alpha_j)$, so that 
\begin{eqnarray*}
T_{t_j}^w(K,\alpha_j)^{t_j} \leq \|P_j\|_{K,w,t_j} \leq \|P_j\|_{K,w,t_j} &\leq& (\|p\|_{K,w,s_{j_0}})^{q_j} \|r_j\|_{K,w,t_j-q_js_{j_0}}  \\  
&\leq&  (L_1-\epsilon)^{q_js_{j_0}}\|r_j\|_K \|w\|_K^{t_j-q_js_{j_0}}.
\end{eqnarray*}
Since $\{\ell_j\}$ is finite, we need only select the $r_j$'s from a fixed finite collection when constructing the 
$P_j$'s; hence we may assume the quantity $\|r_j\|_K$ is uniformly bounded in $j$.  Taking $t_j$-th roots in the above inequality, and the limit as $j\to\infty$,  yields 
$$L_2 = \lim_{j\to\infty} T_{t_j}^w(K,\alpha_j)^{t_j}  \leq L_1-\epsilon.$$  Since $\epsilon>0$ was arbitrary, we are done.
\end{proof}

\begin{definition}\rm
For each $\theta\in\Delta^{\circ}$ we call $T^{w}(K,\theta)$ the \emph{directional Chebyshev constant of $K$} (with weight $w$) associated to $\theta$.  The function $\theta\mapsto\log T^w(K,\theta)$ will be called the \emph{Chebyshev transform of $(K,w)$}.    
\end{definition}

\begin{remark}\rm
The Chebyshev transform takes a set $K$ and weight $w$, and outputs a real-valued function on $\Delta^{\circ}$.  Translated to the complex geometric setting, this corresponds to taking a Hermitian metric on the line bundle and transforming it to a function on the (interior of the) Okounkov body.  Apart from this change in terminology, Nystr\"om's definition in \cite{nystrom:transforming} is the same as above.
\end{remark}

\begin{lemma} \label{lem:310}
\begin{enumerate}
\item The Chebyshev transform $\theta\mapsto \log T^w(K,\theta)$ is convex on $\Delta^{\circ}$. Hence it is either continuous or identically $-\infty$ on $\Delta^{\circ}$. \label{lem:310.1}
\item If the Chebyshev transform is continuous, then for any compact $Q\subset\Delta^{\circ}$, the quantity
$$\sup\bigl\{ |\log T^w_k(K,\alpha(j)) - \log T^w(K,\theta_j)|: \ \alpha(j)/k=:\theta_j\in Q \bigr\}$$ goes to zero  as $k\to\infty$.    \label{lem:310.2}
\end{enumerate}
\end{lemma}

\begin{proof}
To prove the first statement of part (\ref{lem:310.1}), it is sufficient to show that for each $\theta,\phi\in\Delta^{\circ}$ and   $t\in[0,1]$,
\begin{equation}\label{eqn:31a}
T^w(K,t\theta + (1-t)\phi) \leq T^w(K,\theta)^tT^w(K,\phi)^{1-t} .
\end{equation}
Fix $\theta,\phi\in\Delta^{\circ}$ and $t\in[0,1]$.  Choose sequences $\{\alpha_j\}$, $\{\beta_j\}$ in $\NN^m$, and $\{s_j\}$, $\{t_j\}$ in $\NN$, such that $|\alpha_j|,|\beta_j|,s_j,t_j\to\infty$ as $j\to\infty$ and
$$
\tfrac{\alpha_j}{|\alpha_j|}\to \theta,\ \tfrac{\beta_j}{|\beta_j|}\to\phi, \ \tfrac{s_j}{s_j+t_j}\to t.
$$
If we put $\gamma_j:=t_j\alpha_j|\beta_j| + s_j\beta_j|\alpha_j|=: \tilde\alpha_j + \tilde\beta_j$, then by Lemma \ref{lem:subm},
\begin{equation}\label{eqn:32a}
T^w(K,\gamma_j) = T^w(K,\tilde\alpha_j+\tilde\beta_j) \leq T^w(K,\tilde\alpha_j)^{|\tilde\alpha_j|/|\gamma_j|}T^w(K,\tilde\beta_j)^{|\tilde\beta_j|/|\gamma_j|}. 
\end{equation}
We have $\tfrac{\tilde\alpha_j}{|\tilde\alpha_j|}=\tfrac{\alpha_j}{|\alpha_j|}$, $\tfrac{\tilde\beta_j}{|\tilde\beta_j|}=\tfrac{\beta_j}{|\beta_j|}$, and further calculation gives
$$
\tfrac{\gamma_j}{|\gamma_j|}\to t\theta+(1-t)\phi,\quad  \tfrac{|\tilde\alpha_j|}{|\gamma_j|}\to t, \quad \hbox{and } \tfrac{|\tilde\beta_j|}{|\gamma_j|}\to 1-t.
$$
Taking the limit as $j\to\infty$ in (\ref{eqn:32a}) yields (\ref{eqn:31a}), proving the first statement of (\ref{lem:310.1}).  The second statement then follows immediately.

We leave the proof of part (\ref{lem:310.2}) as an exercise in the triangle inequality. (Use Lemma \ref{lem:clim} and the fact that the Chebyshev transform is uniformly continuous on $Q$.)
\end{proof}


\begin{notation}\rm
For a positive integer $k$, let
\begin{itemize}
\item[$M_k$] denote the dimension of $\CC[V]_{\leq k}$ (or number of points in $\calN_k$)
\item[$h_k$] denote the number of points in $\calN_{k}\setminus\calN_{k-1}$ (so $h_k=M_k-M_{k-1}$).
\end{itemize}
\end{notation}

 For each $k\in\NN$, let $\calE_k = \{\be_j\}_{j=1}^{M_k}$ be a basis of polynomials for $\CC[V]_{\leq k}$, with distinct leading terms with respect to their local power series expansions, i.e.,
\begin{equation}  \label{eqn:distinct}
\nu(\be_j) \neq \nu(\be_k) \hbox{ whenever } j\neq k.
\end{equation}
Recall that the bases used in the algorithm to construct the Okounkov body have this property.  
We will also assume that $\calE_1\subset\calE_2\subset\cdots$, and put $\calE:=\bigcup_k\calE_k$, which is a basis for $\CC[V]$ with the above property.  

\begin{definition}\rm
A \emph{monic basis of $\CC[V]$} is a basis $\calE$ that satisfies (\ref{eqn:distinct}) as well as the normalization condition  $\tc(\be_j)=1$ for each $\be_j\in\calE$.
\end{definition}




We now define transfinite diameter. 

\begin{definition} \label{def:vand}  \rm
Let $k$ and  $s$ be positive integers with $s\leq M_k$. Given a finite collection of points $\{\zeta_1,\ldots,\zeta_s\}$, define the \emph{Vandermonde matrix}  
\begin{equation*}
VM_{\calE,k}^w(\zeta_1,\ldots,\zeta_s) \ := \     \begin{bmatrix} \be_1(\zeta_1)_{w,k} &   \cdots & \be_1(\zeta_m)_{w,k} \\
\vdots & \ddots & \vdots \\
\be_s(\zeta_1)_{w,k} & \cdots & \be_s(\zeta_s)_{w,k} \end{bmatrix}  \ = \   [\be_j(\zeta_l)_{w}]_{j,l=1}^s  , 
\end{equation*}
as well as the   \emph{Vandermonde determinant} $$VDM^w_{\calE,k}(\zeta_1,\ldots,\zeta_s):= \det \bigl(VM_{\calE,k}^w(\zeta_1,\ldots,\zeta_s)\bigr).$$ 

For a set $K\subset V$, define
\begin{equation}\label{eqn:Vk}
V_{\calE,k}^w(K,s) := \sup\{|VDM_{\calE,k}^w(\zeta_1,\ldots,\zeta_m)|: \{\zeta_1,\ldots,\zeta_s\}\subset K \}.
\end{equation}
The quantity  
$$d_{\calE,k}^w(K):= \left(V_{\calE,k}^w(K,M_k)\right)^{1/(kM_k)} \hbox { is the \emph{$k$-th order diameter of $K$}},$$  and 
$$d^w(K):=\limsup_{k\to\infty} d_{\calE,k}^w(K)  \hbox{ is the \emph{transfinite diameter of $K$}}.$$  
\end{definition}

\begin{remark}\rm
Since we take the absolute value in (\ref{eqn:Vk}), rearranging the order of the rows in the Vandermonde matrix does not make any difference to the definition of $k$-th order diameters.  We can therefore assume, at each stage $k$, that the rows of $VM^w_{\calE,k}$ have been ordered (and relabelled) so that $\nu(\be_j)\prec_l\nu(\be_k)$ whenever $j<k$.
\end{remark}

\begin{theorem}\label{thm:Ok}
Suppose the Chebyshev transform is an integrable function.  Then the transfinite diameter computed with respect to any monic basis $\calE$ is given by $\displaystyle d^w(K) = \lim_{k\to\infty} d_{\calE,k}^w(K)$, i.e., the limit exists, and
\begin{equation} \label{eqn:transfd}
d^w(K)
= \exp\left(\frac{1}{\vol(\Delta)}\int_{\Delta^{\circ}} \log T^w(K,\theta)d\theta \right), 
\end{equation}
where $d\theta$ denotes the usual volume in $\RR^m$. 
\end{theorem}

In what follows, $\calE=\{\be_j\}_{j=1}^{\infty}$ will be assumed to be a monic basis. 
It will also be convenient to introduce some more notation.
\begin{notation} \rm 
\begin{enumerate}
\item The exponent of $\ttt(\be_j)$ will be denoted by $\alpha(j)$.
\item Given a polynomial $p$, write $o(p)$  to denote an (arbitrary) polynomial $q$ for which $\nu(p)\prec \nu(q)$.
\end{enumerate}
\end{notation}

By definition, $\{\alpha(j):j=1,\ldots,M_k\}  = \calN_k$.

\begin{lemma}\label{lem:<}
Fix $k\in\NN$.  The inequality
\begin{equation*}
V_{\calE,k}^w(K,M_k)  \leq M_k!\prod_{j=1}^{M_k} T^w_k(K,\alpha(j))^k
\end{equation*}
holds.
\end{lemma}

\begin{proof}
For each $j\in\{1,\ldots,M_k\}$,  let 
$$\bt_{j,k}(z) = \be_{j}(z) + o(\be_j)   $$ denote the Chebyshev
polynomial for $(k,\alpha(j))$, i.e., $\|\bt_{j,k}\|_K = T^w_k(K,\alpha(j))^k$.  Now choose a set of $M_k$ points such that $V_{\calE,k}^w(K,M_k) = |VDM_{\calE,k}^w(\zeta_1,\ldots,\zeta_{M_k})|$; then 
$$
|VDM_{\calE,k}^w(\zeta_1,\ldots,\zeta_{M_k})|   
=  \left| \det \begin{bmatrix}
\bt_{1,k}(\zeta_1)_{w,k} & \cdots & \bt_{1,k}(\zeta_{M_k})_{w,k} \\
\vdots & \ddots & \vdots \\
\bt_{M_k-1,k}(\zeta_1)_{w,k} & \cdots & \bt_{M_k-1,k}(\zeta_{M_k})_{w,k} \\
\bt_{M_k,k}(\zeta_1)_{w,k} & \cdots & \bt_{M_k,k}(\zeta_{M_k})_{w,k}  
\end{bmatrix}     \right|
$$
by the invariance of the determinant under row operations.  Expanding the determinant and using the triangle inequality, we obtain 
\begin{eqnarray*} 
V_{\calE,k}^w(K,M_k) &\leq& \sum_{s} |\bt_{1,k}(\zeta_{s(1)})|\cdots|\bt_{M_k,k}(\zeta_{s(M_k)}) |  \\
& \leq & \sum_{s} T_{k}^w(K,\alpha(1))^k\cdots T_k^w(K,\alpha(M_k))^k
\end{eqnarray*}
where the sum is over all permutations $s$ of $\{1,\ldots, M_k\}$, and there are $M_k!$ of these in all.
\end{proof}

\begin{lemma}\label{lem:>}
Fix $k\in\NN$.  Then
$$
V_{\calE,k}^w(K,M_k) \geq \prod_{j=1}^{M_k} T_k^w(K,\alpha(j))^{k}.
$$
\end{lemma}

\begin{proof}
For convenience we will drop the subscripts in polynomial evaluation in what follows, writing e.g. $\be_j(\zeta)=\be_j(\zeta)_{w,k}$. 

 Let us also introduce the following notation: for $j=0,\ldots,M_k-1$, set  
$$
W_j(\zeta_1,\ldots,\zeta_{M_k-j}) := \det 
\begin{bmatrix}
\be_{j+1}(\zeta_1) & \cdots & \be_{j+1}(\zeta_{M_k-j}) \\
\vdots & \ddots & \vdots \\
\be_{M_k}(\zeta_1) & \cdots & \be_{M_k}(\zeta_{M_k-j})
\end{bmatrix} .
$$
In particular $W_0(\zeta_1,\ldots,\zeta_{M_k})=VDM^w_{\calE,k}(\zeta_1,\ldots,\zeta_{M_k})$ and $W_{M_k-1}(\zeta_1)=\be_{M_k}(\zeta_1)$.

We start by observing that for any collection of points $\{\zeta_1,\ldots,\zeta_{M_k}\}\subseteq K$, 
$$
V^w_{\calE,k}(K,M_k) \geq W_0(\zeta_1,\ldots,\zeta_{M_k}).
$$
Using this, we can derive an inequality involving the $W_1$-determinant. To this end, fix   $\{\zeta_1,\ldots,\zeta_{M_k-1}\}\subseteq K$.  Then   $W_0(\zeta_1,\ldots,\zeta_{M_k-1},z)$ (where $z\in K$ is arbitrary) is equal to 
$$
\det
\begin{bmatrix}
\bp_1(\zeta_1) & \cdots & \bp_1(\zeta_{M_k-1}) & \bp_1(z) \\
\vdots & \ddots & \vdots & \vdots \\
\be_{M_k}(\zeta_1) & \cdots & \be_{M_k}(\zeta_{M_k-1}) & \be_{M_k}(z)
\end{bmatrix}, 
$$
where $\bp_1 = \be_1 + o(\be_1)$ and the $o(\be_1)$ expression is a linear combination of the form $\sum_{j=2}^{M_k}\alpha_j\be_j$ obtained via row operations on $VM^w_{\calE,k}(\zeta_1,\ldots,\zeta_{M_k-1},z)$, the Vandermonde matrix.  If the matrix is nonsingular there exist row operations that give a polynomial satisfying  $\bp(\zeta_1)=\cdots=\bp(\zeta_{M_k-1})=0$.  Then 
\begin{eqnarray*}
V_{\calE,k}^w(K,M_k) &\geq &
\left|\det
\begin{bmatrix}
0 & \cdots & 0 & \bp_1(z) \\
\be_2(\zeta_1) & \cdots & \be_{2}(\zeta_{M_k-1}) & \be_2(z) \\
\vdots & \ddots & \vdots & \vdots \\
\be_{M_k}(\zeta_1) & \cdots & \be_{M_k}(\zeta_{M_k-1}) & \be_{M_k}(z)
\end{bmatrix} \right| \\
&=&  |\bp_1(z)|\cdot|W_1(\zeta_1,\ldots,\zeta_{M_k-1})|.
\end{eqnarray*}
Choosing $\eta\in K$ such that $|\bp_1(\eta)|=\|\bp_1\|_K\geq T^w_k(K,\alpha(1))^k$, we have
\begin{equation}\label{eqn:W}
V_{\calE,k}^w(K,M_k)\geq  |W_0(\zeta_1,\ldots,\zeta_{M_k-1},\eta)| =    T^w_k(K,\alpha(1))^k |W_1(\zeta_1,\ldots,\zeta_{M_k-1})|.
\end{equation}

Now consider fixing points $\zeta_1,\ldots,\zeta_{M_k-2}$ and carrying out a similar argument as above: construct a polynomial $\bp_2= \be_2 + o(\be_2)$ with $\bp_2(\zeta_1)=\cdots=\bp_2(\zeta_{M_k-2})=0$ using row operations on the matrix for $W_1$, then choose $\eta\in K$ for which $|\bp_2(\eta)|=\|\bp_2\|_K$. This gives the inequality 
\begin{equation}\label{eqn:W1}
|W_1(\zeta_1,\ldots,\zeta_{M_k-2},\eta)|   \geq T_k^w(K,\alpha(2))^k|W_2(\zeta_1,\ldots,\zeta_{M_k-2})|.
\end{equation}
We use (\ref{eqn:W}) to estimate the left-hand side of (\ref{eqn:W1}), observing that the upper bound on the left-hand side of (\ref{eqn:W}) is valid for an arbitrary collection of $M_{k-1}$ points of $K$.  Hence 
$$V_{\calE,k}^w(K,M_k) \geq T_k^w(\alpha(1))^k T_k^w(K,\alpha(2))^kW_2(\zeta_1,\ldots,\zeta_{M_k-2}) .$$
Now it is easy to see that the argument can be iterated to obtain the estimate 
$$V_{\calE,k}^w(K,M_k) \geq \Bigl(\prod_{j=1}^s T_k^w(K,\alpha(j))^k \Bigr) W_s(\zeta_1,\ldots,\zeta_{M_k-s})$$ for successive values  $s=3,4,\ldots$.  The lemma is proved when $s=M_k-1$.
\end{proof}

Theorem \ref{thm:Ok} is proved by putting together the above lemmas.

\begin{proof}[Proof of Theorem \ref{thm:Ok}]  We have by Lemmas \ref{lem:<} and \ref{lem:>}, 
$$ 
 \prod_{j=1}^{M_k} T_k^w(K,\alpha(j))^{k}    \leq	  V_{\calE,k}^w(K,M_k)   \leq  M_k!\prod_{j=1}^{M_k} T^w_k(K,\alpha(j))^k .
$$
Let us take $kM_k$-th roots in this inequality and let $k\to\infty$.  Since   
$(M_k!)^{1/(kM_k)}\to 1$ as  $k\to\infty$,  
$$ 
d^w(K) = \lim_{k\to\infty} V_{\calE,k}^w(K,M_k)^{1/(kM_k)} =\lim_{k\to\infty} \left(\prod_{j=1}^{M_k} T^w_k(K,\alpha(j))^k\right)^{1/(kM_k)}.  $$
We need to show that the limit on the right-hand side of the above converges to the right-hand side of (\ref{eqn:transfd}).  To see this, observe that the limit may be rewritten as
\begin{equation}\label{eqn:pf312}
\lim_{k\to\infty} \exp\biggl(   
\frac{1}{M_k}\sum_{j=1}^{M_k} \log T^w_k(K,\alpha(j)) 
\biggr).
\end{equation}
We look at the limit of the expression inside the parentheses.  If we fix a compact convex body $Q\subset\Delta^{\circ}$, then by Lemma 3.10(2), we have 
$$
\lim_{k\to\infty} \frac{1}{M_k}\sum_{\alpha(j)\in kQ}^{M_k} \log T^w_k(K,\alpha(j)) \ = \ \lim_{k\to\infty} \frac{1}{M_k}\sum_{\alpha(j)\in kQ} \log T^w(K,\theta_j),
$$ 
where $\theta_j=\alpha(j)/k$.  Now (2.12) in Proposition 2.15 implies that the discrete measure $\frac{1}{M_k}\sum_{j=1}^{M_k}\delta_{\theta_j}$ converges weak-$^*$ on $\Delta^{\circ}$ to the uniform measure  $\frac{1}{\vol(\Delta)}\,d\theta$ as $k\to\infty$.\footnote{(2.12) holds for compact convex bodies, and these generate all Borel sets.  The claim then follows by standard measure theory.}    Hence 
$$
\lim_{k\to\infty} \frac{1}{M_k}\sum_{\alpha(j)\in kQ} \log T^w(K,\theta_j) \  = \   \frac{1}{\vol(\Delta)}\int_Q\log T^w(K,\theta)\, d\theta.
$$
Since $Q$ was arbitrary, one can consider a sequence of compact convex sets increasing to $\Delta^{\circ}$ and obtain
$$
\lim_{k\to\infty} \frac{1}{M_k}\sum_{j=1}^{M_k} \log T^w(K,\theta_j) \  = \   \frac{1}{\vol(\Delta)}\int_{\Delta^{\circ}}\log T^w(K,\theta)\, d\theta,
$$
using a standard convergence argument.
\end{proof}

\begin{remark}\rm
If $T^w(K,\theta)$ is nonzero and uniformly bounded from above then the Chebyshev transform is integrable, since the Okounkov body $\Delta$ is compact.  It is straightforward to get an upper bound when $K$ is compact, or when $w$ is admissible.  
\end{remark}

Nystr\"{o}m's formula in \cite{nystrom:transforming} is actually stated in terms of a ratio.  
The analogous result in our setting is to take another set $L$ and admissible weight $u$, and compute the associated quantities.  It is easy to see that  
$$
\frac{d^w(K)}{d^u(L)} \ = \ 
\exp\left( \frac{1}{\vol(\Delta)}\int_{\Delta^{\circ}} (\log T^w(K,\theta) - \log T^u(L,\theta)) d\theta \right). 
$$
In place of the ratio of transfinite diameters on the left-hand side, the formula in \cite{nystrom:transforming} has a \emph{mixed Monge-Amp\`ere energy}, which for the above case is
$$\calE(P_K(\varphi),P_L(\psi)) := \int_{V} (P_K(\varphi)-P_L(\psi))\sum_{j=0}^m (dd^cP_K(\varphi))^j\wedge (dd^cP_L(\psi))^{m-j}.$$
Here, $P_K(\varphi),P_L(\psi)$ are plurisubharmonic extremal functions defined in terms of the sets $K,L$ and functions $\varphi=-\log w,\psi=-\log u$.\footnote{From the point of view of pluripotential theory, it is more natural to consider the logs of the weights.} 

The Monge-Amp\`ere energy appears in Nystr\"{o}m's formula because his proof is based on pluripotential theory and the $L^2$ theory associated to Bernstein-Markov measures.  It is equivalent to the above formula because of the (highly nontrivial) equality $$\frac{d^w(K)}{d^u(L)}=\calE(P_K(w),P_{L}(u)),$$ which is a version of \emph{Rumely's formula} proved in \cite{bermanboucksom:growth} (see also \cite{levenberg:bermanboucksom} for a self-contained exposition). The original statement and proof of Rumely's formula is in \cite{rumely:robin}.

\section{Comparison with classical results} \label{sec:comparison}


The arguments in the previous section are almost the same as those of Zaharjuta in \cite{zaharjuta:transfinite}.  He derived an integral formula for the classical transfinite diameter of a compact set $K\subset\CC^n$, denoted $d(K)$.  Let $\calE=\{z^{\alpha(j)}\}_{j=1}^{\infty}$ be the enumeration of the monomials in $z_1,\ldots,z_n$ according to the \emph{grevlex order} on multi-indices $\alpha(j)=(\alpha(j)_1,\ldots,\alpha(j)_n)$, which is defined as follows: if $j<k$ then 
\begin{itemize}
\item either $|\alpha(j)|<|\alpha(k)|$; or 
\item $|\alpha(j)=|\alpha(k)|$, and there exists $\nu\in\{1,\ldots,n\}$ such that 
$$\alpha(j)_{1}=\alpha(k)_{1},\ldots,\alpha(j)_{\nu-1}=\alpha(k)_{\nu-1}, \quad \hbox{and } \alpha(j)_{\nu}<\alpha(k)_{\nu}. $$
\end{itemize}
 

 Let  
 \begin{equation} \label{eqn:simplex}
 \Sigma = \{\theta=(\theta_1,\ldots,\theta_n)\in\RR^n: \ \theta_k\geq 0\ \forall\, k,\  \sum_k\theta_k = 1\} \end{equation}
  be the standard $(n-1)$-dimensional simplex in $\RR^n$, and $\Sigma^{\circ} := \{\theta \in\Sigma: \theta_k>0\ \forall\,k\}$ its interior. Zaharjuta showed that Chebyshev constants parametrized by $\theta\in\Sigma^{\circ}$ can be defined as follows: let 
\begin{equation}\label{eqn:tauK}
\tau_j(K):=\inf\Bigl\{\|p\|_K^{1/|\alpha(j)|}: p(z)=z^{\alpha(j)} + \sum_{\nu=1}^{j-1}c_{\nu}z^{\alpha(\nu)} \Bigr\}, \ \hbox{and} \ \tau(K,\theta) := \lim_{\substack{j\to\infty\\ \frac{\alpha(j)}{|\alpha(j)|}\to\theta }} \tau_j(K).  
\end{equation}
\begin{remark} \label{rmk:tau} \rm 
We will also use the notation 
$\tau_{\alpha}(K) = \tau_{\alpha(j)}(K) = \tau_j(K)$
when $\alpha=\alpha(j)$ for some positive integer $j$.
\end{remark}

 Next consider $V_{\calE,k}(K,M_k)$ as in Definition \ref{def:vand} (we suppress the dependence on the trivial weight $w\equiv 1$), and let
 \begin{equation}\label{eqn:lk}
 L_k:=\sum_{s=1}^k sh_s.
 \end{equation}
    Zaharjuta showed that the limit $$\displaystyle d(K):=\lim_{k\to\infty} \bigl(V_{\calE,k}(K,M_k) \bigr)^{1/L_k},$$ the \emph{transfinite diameter of $K$}, exists and satisfies the following formula.  
\begin{theorem} \label{thm:z}
We have $\displaystyle d(K) = \exp\left( \frac{1}{\vol(\Sigma)}\int_{\Sigma^{\circ}} \log\tau(K,\theta)d\theta   \right)$.
\end{theorem}

Later, Jedrzejowski \cite{jedrzejowski:homogeneous} showed a similar formula for the \emph{homogeneous} transfinite diameter  in $\CC^n$; let us denote it by $\hat d(K)$. The homogeneous transfinite diameter  for a compact set $K$ is defined by  
$$\hat d(K):=\lim_{k\to\infty} \left(\sup\bigl\{|VDMH(\zeta_1,\ldots,\zeta_{h_k})|^{1/(kh_k)}: \{\zeta_1,\ldots,\zeta_{h_k} \}\subset K \bigr\} \right),  $$ where, with $\be_j:=z^{\alpha(j)}$,  
$$
VDMH(\zeta_1,\ldots,\zeta_{h_k}) = \det\begin{bmatrix}
\be_{M_{k-1}+1}(\zeta_1) &  \cdots & \be_{M_{k-1}+1}(\zeta_{h_k})  \\
\vdots & \ddots & \vdots \\
\be_{M_k}(\zeta_1) & \cdots & \be_{M_k}(\zeta_{h_k})   
\end{bmatrix}.
$$
One can construct homogeneous Chebyshev constants $\hat\tau(K,\theta)$ as limits of constants $\hat\tau_j(K)$, where the latter are defined as in (\ref{eqn:tauK}), but with the inf restricted to homogeneous polynomials  (i.e.  $c_{\nu}=0$ if $|\alpha(\nu)|<|\alpha(j)|$). The homogeneous formula is the following.
\begin{theorem}\label{thm:j} 
We have $\displaystyle \hat d(K) = \exp\left( \frac{1}{\vol(\Sigma)}\int_{\Sigma^{\circ}} \log\hat\tau(K,\theta)d\theta   \right)$.
\end{theorem}

\begin{remark} \rm
Bloom and Levenberg have studied weighted versions of directional Chebyshev constants \cite{bloomlev:weighted}. For a fixed admissible weight $w$, $\tau_j^w(K)$ is defined as in (\ref{eqn:tauK}) by a sup over quantities of the form  $\|w^{|\alpha(j)|}p\|_K$.  A weighted generalization of Theorem \ref{thm:z} was proved later in \cite{bloomlev:transfinite}.
\end{remark}

In what follows, we expand a bit on the relationship between the above theorems and Theorem \ref{thm:Ok}.  

\subsection{Homogeneous transfinite diameter}  
Consider $\CC^n$ with coordinates given by $(z_0,z_1,\ldots,z_{n-1})$, and consider a compact subset of the form $\{1\}\times K$ (i.e. $z_0=1$), where $K$ is compact in $\CC^{n-1}$.  We will also assume for what follows that $K$ avoids the hyperplane $z_{n-1}=0$. We describe the relation between Theorems \ref{thm:Ok} and \ref{thm:j}.  

Consider the monomials in $\CC[z_0,\ldots,z_{n-1}]$ with the ordering $\prec$ defined as the lexicographic order for which $z_{n-1}\prec\cdots\prec z_{1}\prec z_0$.  The homogeneous polynomials $\{\be_j\}_{j=M_{k-1}+1}^{M_k}$ of degree $k$ are given by  $$\be_{M_{k-1}+1}=z_{n-1}^k,\ \be_{M_{k-1}+2}=z_{n-2}z_{n-1}^{k-1},\ \ldots, \ \be_{M_{k}}=z_0^{k}.$$ 
On the variety $V:=\{z_0=1\}$, polynomials are given by $\CC[z_1,\ldots,z_{n-1}]$ and for $z=(1,z_1,\ldots,z_{n-1})\in \{1\}\times K\subset V$ we can identify $\CC[z_0,\ldots,z_{n-1}]_k$ (the homogeneous polynomials of degree $k$ in $n$ variables) with $\CC[V]_{\leq k}\simeq\CC[z_1,\ldots,z_{n-1}]_{\leq k}$, via 
$$
z_0^{\alpha_0}z_1^{\alpha_1}\cdots z_{n-1}^{\alpha_{n-1}} = z_1^{\alpha_1}\cdots z_{n-1}^{\alpha_{n-1}} \quad 
(  \alpha_0+\cdots+\alpha_{n-1}=k)   .
$$
 We also have the relation to $\CC[v_0,v_1,\ldots,v_{n-2}]$ given by
$$
|z_1^{\alpha_1}\cdots z_{n-1}^{\alpha_{n-1}}| = |w(v)^kv_0^{\alpha_0}v_1^{\alpha_1}\cdots v_{n-2}^{\alpha_{n-2}}|
$$
where $v$-coordinates are given by 
\begin{equation}\label{eqn:vc} v_0 =1/z_{n-1}, \quad v_j=z_j/z_{n-1} \ \hbox{ for }j=1,\ldots,n-2,\end{equation} and weight $w$ defined by the formula $w(v)= 1/z_{n-1}$.      Identifying $z_1^{\alpha_1}\cdots z_{n-1}^{\alpha_{n-1}}\in\CC[V]_{\leq k}$ with $v_0^{\alpha_{n-1}}v_1^{\alpha_2}\cdots v_{n-2}^{\alpha_{n-2}}\in\CC[v_0,v_2,\ldots,v_{n-1}]_{\leq k}$, polynomial evaluation is related by $|p(z)| = |w(v)^kp(v)|$.

\smallskip

Now observe that
\begin{itemize}
\item[(i)] $h_k$ (computed with respect to $\CC[z_0,z_1,\ldots,z_n]$) is the same as $M_k$ for the spaces 
$\CC[V]_{\leq k}$ and $\CC[v_0,v_1,\ldots,v_{n-2}]_{\leq k}$.
 \item[(ii)] Under the identifications described above, the lex order $\prec_l$ on $\CC[z_0,\ldots,z_n]_k$ translates to lex order on the corresponding monomials in $\CC[v_0,v_1,v_{n-2}]_{\leq k}$, and to the \emph{reverse} of grevlex order (which we will denote by $\prec_g$) on the corresponding monomials in affine coordinates on $V=\{z_0=1\}$. An example of a pair of monomials when $n=4$ and $k=5$, together with the corresponding pairs in the other spaces, is 
 \begin{eqnarray*}
 z_0^3z_1z_3 \succ_l z_0^2z_1^2z_4 \ \hbox{ in }\CC[z_0,z_1,z_2,z_3,z_4]  &\sim&  v_0^3v_1v_3\succ_l v_0^2v_1^3 \ \hbox{ in }\CC[v_0,v_1,v_2,v_3]  \\
  &\sim&  z_1z_3\prec_{g} z_1^2z_4  \  \hbox{ in } \CC[z_1,z_2,z_3,z_4].
 \end{eqnarray*}
 \end{itemize}
   These identifications yield, for any collection of points $\{\zeta_1,\ldots,\zeta_{h_k}\}\subset V$,   
$$VDMH(\zeta_1,\ldots,\zeta_{h_k}) = VDM_{\calE,k}^w(\zeta_1,\ldots,\zeta_{h_k})$$ 
where the right-hand side may be interpreted either in $z$-coordinates with $w\equiv 1$ (unweighted) and $\calE$ being monomials in $z$, or in $v$-coordinates with monomials in $v$ and the weight $w=w(v)$ described above.  The equivalence of the above determinants, and the fact that the same root is taken at each stage ($kh_k$-th or $kM_k$-th, see (i) above), yields the equality $\hat d(\{1\}\times K) = d^w(\{1\}\times K)$.
 
 Chebyshev constants are also related.  (For convenience of notation, let us simply write $K$ for $\{1\}\times K$ in what follows.)  Given $\theta\in\Sigma$, i.e.,  $\theta=(\theta_0,\ldots,\theta_{n-1})$ with $\theta_j>0$ and $\sum\theta_j =1$, we have 
 \begin{equation}\label{eqn:tt}
 \hat\tau(K,\theta) = T^w(K,\theta')
 \end{equation}
 where we interpret the right-hand side in $v$ coordinates, in which $w=w(v)$ and $\theta'=(\theta_0,\ldots,\theta_{n-2})$.\footnote{Since the Okounkov body is constructed with reference to local coordinates, Theorem \ref{thm:Ok} only applies directly to the $v$-coordinate setting.}  
As $\theta$ varies over $\Sigma^{\circ}$, $\theta'$ varies over the projection of these points to the first $n-1$ coordinates, which fills the interior of the region $S$ in $\RR^{n-1}$ given by   $\theta_0+\cdots+\theta_{n-2}\leq 1$, $\theta_j\geq 0$.  Now $d\theta'$, the volume in $\RR^{n-1}$, is the push forward of $d\theta$ (the $n-1$-dimensional volume in the plane containing $\Sigma$), scaled by a factor of $\vol(\Sigma)/\vol(S)$. Hence 
 $$ \frac{1}{\vol(\Sigma)}\int_{\Sigma^{\circ}}\hat\tau(K,\theta)d\theta \ = \ \frac{1}{\vol(S)}\int_{S^{\circ}}T^w(K,\theta')d\theta'.$$  
Since $S$ is the Okounkov body of $\CC^{n-1}$, Theorem \ref{thm:j} is just Theorem \ref{thm:Ok} under a change of variable.

 \subsection{Transfinite diameter}
 Theorem \ref{thm:z} is slightly different, but closely related.  We describe its precise relationship to Theorem \ref{thm:Ok} in what follows.  We will work with a compact set $K\subset\CC^{n-1}$ (with variables $z_1,\ldots,z_{n-1}$) and reuse the material above, relating $K$ to $\{1\}\times K\subset\CC^n$.  In particular, $z_0$ denotes the additional variable, and $v$-coordinates are defined as in (\ref{eqn:vc}).  It is easy to see that $d(K)$ is related to $\hat d(\{1\}\times K)$: we have  
$$
VDMH(\zeta_1,\ldots,\zeta_{M_k}) = VDM_{\calE,k}(\zeta_1',\ldots,\zeta_{M_k}')
$$
where $M_k$ is counted with respect to $\CC[z_1,\ldots,z_{n-1}]$ (and therefore coincides with `\emph{$h_k$ for $\CC[z_0,\ldots,z_{n-1}]$}'), and we use the notation $\CC^n\ni\zeta_j=(1,\zeta_j')\in\CC\times\CC^{n-1}$.  

We saw previously that $\hat d(\{1\}\times K) = d^w(K)$, but  $d(K)$ in Theorem \ref{thm:z} takes a slightly different root at each step $k$, which affects the normalization.  We omit the calculation of the relevant limit.
\begin{proposition}\label{prop:46}
We have $\displaystyle \lim_{k\to\infty} \, \frac{L_k}{kM_k} \  =  \  \frac{n-1}{n}$, where $L_k$ is as in (\ref{eqn:lk}).  Hence
$$
d^w(K) =  \hat d(\{1\}\times K) =  d(K)^{(n-1)/n},
$$
where the above quantities are from Theorems \ref{thm:Ok}, \ref{thm:j} and \ref{thm:z} respectively. \qed 
\end{proposition}

We now turn to the integral formulas.  Suppose $\theta\in S^{\circ}$, and consider the constants $T^w(K,\theta)$ defined as on the right-hand side of (\ref{eqn:tt}).  These constants (associated to  $v$-coordinates)  transform to constants $T(K,\tilde\theta)$ associated to $z$-coordinates, with parameters related by
\begin{equation}\label{eqn:tildetheta}
S^{\circ}\ni\theta = (\theta_0,\ldots,\theta_{n-2}) \iff  (\theta_1,\ldots,\theta_{n-1}) = \tilde\theta \in \tilde S^{\circ},
\end{equation}
whenever $\theta_0+\cdots+\theta_{n-1}=1$.\footnote{We will write $\tilde S$ for the standard triangle associated to $z$-coordinates (and write  $\tilde S^{\circ}$ for its interior) but use the same labels $K,p$ for sets and polynomials.} 
 The constants $T(K,\tilde\theta)$  satisfy the following homogeneity property.
 \begin{lemma}
Let $\theta\in \tilde S^{\circ}$.  Then for all $t\in[0,1]$,
$$
\log T(K,t\theta) = t\log T(K,\theta).
$$
Moreover, if $r\theta=\varphi\in\Sigma^{\circ}$ for some $r>1$ then
$$
r\log T(K,\theta) =\log\tau(K,\varphi)
$$
where $\tau(K,\varphi)$ is the directional Chebyshev constant given by (\ref{eqn:tauK}).
\end{lemma}

\begin{proof}
Fix a positive integer $k$, and a positive integer $j<M_k$.  Define the Okounkov body with respect to $v$-coordinates as above, and let $w$ be the corresponding weight.  Let $p$ be a Chebyshev polynomial such that $T_k^w(K,\alpha(j))^k=\|w^kp\|_K$.  In $z$-coordinates, the weight becomes trivial and the sup translates to a sup over quantities of the form $\|p\|_K$,
which we will denote by $T_k(K,\tilde\alpha(j))^k$.  (The exponent $\tilde\alpha(j)$ is such that $z^{\tilde\alpha(j)}$ corresponds to $v^{\alpha(j)}$.)  It is easy to see that $p$ must be a Chebyshev polynomial of degree $s=|\tilde\alpha(j)|$ (i.e., one that attains the inf in  (\ref{eqn:tauK})  for $\tau_j(K)=\tau_{\tilde\alpha(j)}(K)$),  and so
\begin{equation} \label{eqn:tauT}
\tau_j(K)^{s} = T_k(K,\tilde\alpha(j))^k.
\end{equation}
Suppose, for some rational number $\frac{a_j}{b_j}\in(0,1)$, that the $(n-1)$-tuple $\frac{a_j}{b_j}\tilde\alpha(j)$ has integer entries (i.e. $b_j$ divides all components).  Then $\frac{a_j}{b_j}\tilde\alpha(j)= \tilde\alpha(j')$ for some $j'<j$.  It is easy to see that $|\tilde\alpha(j)|=s$ implies $|\tilde\alpha(j')|=\frac{a_j}{b_j}s$.  Similar to (\ref{eqn:tauT}), one also has 
\begin{equation}\label{eqn:tauT1}
\tau_{j'}(K)^{\frac{a_j}{b_j}s} = T_k(K,\tfrac{a_j}{b_j}\tilde\alpha(j))^k.
\end{equation}
Equations (\ref{eqn:tauT}) and (\ref{eqn:tauT1}) yield the lemma after some analysis.  Precisely, take sequences of positive integers $a_j,b_j$ such that $b_j\to\infty$ and $a_j/b_j\to t\in(0,1)$. Then take a sequence of exponents $\alpha_j$ such that $b_j$ divides each component of $\alpha_j$, and a sequence of integers $k_j>|\alpha_j|=:s_j$ such that $\alpha_j/k_j\to\theta\in S^{\circ}$.  Let $\beta_j:=\frac{a_j}{b_j}\alpha_j$; then $\beta_j/k_j\to t\theta$ where $\beta_j=\frac{a_j}{b_j}\alpha_j$.   

Define $\varphi:=\lim_{j\to\infty} \alpha_j/|\alpha_j|$.  Then $\varphi=r\theta\in\Sigma$, where $r=\lim_{j\to\infty} k_j/s_j$.   For the same reason, $\varphi=\lim_{j\to\infty}\beta_j/|\beta_j|$ also holds, and as $j\to\infty$, 
\begin{eqnarray}
 \label{eqn:46}  T_{k_j}(K,\beta(j))^{{(b_jk_j)}/{(a_js_j)}} &=& \tau_{\beta_j}(K)  \to \tau(K,\varphi), \hbox{ and} \\
  \label{eqn:47} T_{k_j}(K,\alpha(j))^{k_j/s_j} &=& \tau_{\alpha_j}(K) \to \tau(K,\varphi).
\end{eqnarray}
Now note that as $j\to\infty$ we also have $T_{k_j}^w(K,\beta(j))^{{(b_jk_j)}/{(a_js_j)}}\to T^w(K,t\theta)^{r/t}$ and $T_{k_j}^w(K,\alpha(j))^{k_j/s_j}\to T(K,\theta)^r$ on the left-hand sides.  By (\ref{eqn:46}) and (\ref{eqn:47}) we may equate all of the asymptotic quantities and take logs.  Both statements of the lemma follow immediately. 
\end{proof}

Using the above lemma, we can directly relate the integrals of Theorems \ref{thm:Ok} and \ref{thm:z}.  First, we transform the integral in Theorem \ref{thm:Ok} to $z$-coordinates: 
$$
\frac{1}{\vol(S)}\int_{S^{\circ}}\log T^w(K,\theta)\, d\theta  \ = \ \frac{1}{\vol(\tilde S)}\int_{\tilde S^{\circ}}\log T(K,\tilde\theta)\,d\tilde\theta .
$$
Now observe that $\tilde S$ can be expressed as the union $\tilde S=\bigcup_{t\in[0,1]} t\Sigma$ where $\Sigma\in\RR^{n-1}$ is the $(n-2)$-dimensional simplex defined as in (\ref{eqn:simplex}).   
Using the map $$(0,1)\times\Sigma^{\circ}\ni(t,\theta')\mapsto t\theta'=\tilde\theta\in \tilde S^{\circ},$$ we see that the volume element $d\tilde\theta$ on $\tilde S^{\circ}$ may be decomposed as $t^{n-2}{d\theta'}dt$.  Continuing, we have
\begin{eqnarray}
\frac{1}{\vol(\tilde S)}\int_{\tilde S^{\circ}} \log T(K,\tilde\theta)\,d\tilde\theta  &=& \frac{1}{\vol(\tilde S)}\int_0^1\int_{\Sigma^{\circ}} \log T(K,t\theta')t^{n-2}d\theta'dt  \nonumber  \\
&=& \frac{1}{\vol(\tilde S)}\int_0^1 \int_{\Sigma^{\circ}} \log \tau(K, \theta')t^{n-1}d\theta '    \nonumber  \\
&=& \frac{1}{\vol(\tilde S)} \int_0^1t^{n-1}dt   \int_{\Sigma^{\circ}}\log\tau(K,\theta')d\theta'   \nonumber  \\ 
& = & \frac{1}{n\vol(\tilde S)} \int_{\Sigma^{\circ}}\log\tau(K,\theta')d\theta', \label{eqn:48}
\end{eqnarray}
where we use the previous lemma.  We compute $\vol(\tilde S)$ using the same decomposition:
$$
\vol(\tilde S) = \int_0^1 \vol(t\Sigma)dt = \int_0^1 t^{n-2}\vol(\Sigma)dt = \frac{1}{n-1}\vol(\Sigma).
$$ 
Finally, substitute the above expression for $\vol(\tilde S)$ into (\ref{eqn:48}).  Altogether, we have
$$
\frac{1}{\vol(S)}\int_{S^{\circ}}T^w(K,\theta)d\theta \ = \ \left(\frac{n-1}{n}\right)\frac{1}{\vol(\Sigma)}\int_{\Sigma^{\circ}}\tau(K,\theta')d\theta'.
$$
Observe that the normalization agrees with Proposition \ref{prop:46}.

\section{Further properties}  \label{sec:further}

In this section, we study further properties of Chebyshev constants and associated notions.  Specific results will be given on the complexified sphere.  More general results will be the subject of future research.

\subsection{General collections of polynomials.} 
We want to reuse the methods of Theorem \ref{thm:Ok} in a more general context, so let us extract the essential ingredients required for the proof.  The Vandermonde determinant used in the limiting process is defined in terms of a collection of polynomials $\calE=\{\be_j\}_{j=1}^{\infty}$ with some additional structure related to a grading with respect to multiplication: $\calE=\bigcup_{k=0}^{\infty}\calE_k$ with $\calE_0\subset\calE_1\subset\cdots$, and $\calE_k\calE_l = \{pq: p\in\calE_k, q\in\calE_l\}$ is contained in the span of $\calE_{k+l}$ for each pair of non-negative integers $k,l$.  The structure associated to this grading allows us to compute the Okounkov body and associated Chebyshev constants, and consists of the following things.
\begin{enumerate}
\item There is a function $\nu:\span(\calE)\to\NN_0^m$ which is one-to-one on $\calE$, and an associated convex body   given as follows. Let $\calN_k=\{\nu(\be_j)\}_{j=1}^{M_k}$ where $\calE_k=\{\be_1,\ldots,\be_{M_k}\}$. Then let $\Delta_{\calE,k}$ be the convex hull of $\frac{1}{k}\calN_k$ and let $\Delta_{\calE}$ be the convex hull of $\bigcup_k\Delta_{\calE,k}$.  
\item For each $k\in\NN$ and $\alpha\in\calN_k$ there is a class of polynomials $\calM(k,\alpha)\subset\span(\calE_k)$ such that $\be_j\in\calM(k,\alpha)$ if $\nu(\be_j)=\alpha$.  These classes satisfy the properties 
\begin{eqnarray*}
 & & \calM(k,\alpha) + p \subseteq \calM(k,\alpha)\quad \hbox{if } p\in\span(\calE_k) \hbox{ and } \nu(p)\prec\alpha, \\
 & & \calM(k_1,\alpha_1)\calM(k_2,\alpha_2) \subseteq \calM(k_1+k_2,\alpha_1+\alpha_2).  
\end{eqnarray*}
\item The discrete measure $\frac{1}{M_k}\sum_{\alpha\in\calN_k} \delta_{\alpha/k}$ converges weak-$^*$ to $\frac{1}{\vol(\Delta_{\calE})}d\theta$, where $d\theta$ is the usual volume measure on $\RR^m$ restricted to $\Delta_{\calE}$.
\end{enumerate}
In addition, there is a weighted polynomial evaluation with respect to some admissible weight $w:K\to\CC$, with  
$$
p(z)_{w,k}q(z)_{w,l} = (pq)(z)_{w,k+l}.
$$

With these properties, one can then construct, for $\alpha\in\calN_k$, 
 $$ T_k(K,\alpha) :=  \inf\{ \|p\|_{K,w,k}: p\in\calM(\alpha,k)\}^{1/k} $$
(we suppress the dependence on $w$ here and in what follows), as well as the function on the interior of the convex body,  $T:\Delta_{\calE}^{\circ}\to[0,\infty)$ given by 
$$
T(K,\theta) := \lim_{ \substack{k\to\infty \\ \alpha/k\to\theta}} T_k(K,\alpha).
$$
Then defining   
$$d_{\calE,k}(K):=\sup\bigl\{|VDM_{\calE,k}(\zeta_1,\ldots,\zeta_{M_k})|^{1/(kM_k)}:\  \{\zeta_1,\ldots,\zeta_{M_k}\}\subset K\bigr\}$$
as in Definition \ref{def:vand}, we have 
\begin{equation} \label{eqn:Okgen}
d_{\calE}(K) = \exp\left( \frac{1}{\vol(\Delta_{\calE})} \int_{\Delta_{\calE}^{\circ}} \log T(K,\theta)d\theta \right). 
\end{equation}

\begin{remark}\rm 
  The main point is that such a formula arises for \emph{any} collection $\calE$ of polynomials with the type of structure given above; for example, the basis of a graded subalgebra of $\CC[V]$.   (In the complex geometric setting, Okounkov bodies associated to subalgebras of global sections have been studied e.g., in \cite{hisamoto:on}.) 
\end{remark}

\subsection{Monomials on the sphere}\label{sec:52}
Consider again the complexified sphere $V\subset\CC^3$ given by (3.1)
  which is spanned by monomials of the form $z_1^{\alpha_1}z_2^{\alpha_2}z_3^{\alpha_3}$ with $\alpha_1,\alpha_2\in\NN_0$ and $\alpha_3\in\{0,1\}$.  Let $\calE$ denotes the subcollection of monomials of the form $z_1^{\alpha_1}z_2^{\alpha_2}$ that span $\CC[z_1,z_2]\subseteq\CC[V]$; as usual, for $\alpha=(\alpha_1,\alpha_2)$, write $z^{\alpha}:=z_1^{\alpha_1}z_2^{\alpha_2}$, and order the monomials by grevlex.  We treat them as in the classical theory:
\begin{enumerate}
\item let $M_k$ denote the number of monomials in $\calE$ of degree $\leq k$;
\item let $\nu:\CC[z_1,z_2] \to \NN_0^2$ return the leading exponent of a polynomial (with respect to grevlex); 
\item let $\calM(\alpha,k)\subset\CC[z_1,z_2]_{\leq k}$ denote the class of monic polynomials of the form $$z^{\alpha} + (\hbox{lower terms wrt grevlex});$$ 
\item and let $VDM(\zeta_1,\ldots,\zeta_{M_k})$ denote the Vandermonde determinant associated to $\calE$ with the standard polynomial evaluation at affine points.
\end{enumerate}
Define Chebyshev constants and transfinite diameter:
\begin{eqnarray*}  T_{\calE,k}(K,\alpha) &:=& \inf\{\|p\|_K:\ p\in\calM(k,\alpha)   \}^{1/k}, \  \hbox{and}  \\  
d_{\calE,k}(K) &:=& \sup\{|VDM(\zeta_1,\ldots,\zeta_{M_k})|^{1/(kM_k)}: \zeta_j\in K \ \forall\, j\}.
\end{eqnarray*}
The Okounkov body $\Delta_{\calE}$ associated to $\calE$ is the standard triangle,
 $$\Delta_{\calE} = S=\{\theta=(\theta_1,\theta_2)\in\RR^2: \theta_1,\theta_2\geq 0, \theta_1+\theta_2\leq 1\}.$$
The limits
$$
d_{\calE}(K):= \lim_{k\to\infty} d_{\calE,k}(K) \  \hbox{ and } \  
  T_{\calE}(K,\theta):= \lim_{\substack{ k\to\infty\\ \alpha/k\to\theta  }} T_{\calE,k}(K,\alpha)
$$
exist (for all $\theta\in S^{\circ}$ in the latter), and
\begin{equation}\label{eqn:thm51}
d_{\calE}(K) = \exp\left( \frac{1}{\vol(S)}\int_{S^{\circ}} \log T_{\calE}(K,\theta)\,d\theta  \right) .   
\end{equation} 

Observe that $\calE$ (and hence $\CC[z_1,z_2]$) cannot distinguish between points with different $z_3$-coordinates: if $\pi:V\to\CC^2$ is the projection $(z_1,z_2,z_3)\mapsto(z_1,z_2)$,  then for all $\be_j\in\calE$, we have $\be_j(\zeta)= \be_j(\eta)$ whenever $\pi(\zeta)=\pi(\eta)$.  As an immediate consequence, for any collections $\{\zeta_1,\ldots,\zeta_{M_k}\}$ and $\{\eta_1,\ldots,\eta_{M_k}\}$ of points of $V$, 
\begin{equation}\label{eqn:vdme}
VDM(\zeta_{1},\ldots,\zeta_{M_k}) = VDM(\eta_1,\ldots,\eta_{M_k}) \hbox{ whenever }  \eta_j\in\pi^{-1}(\{\pi(\zeta_j)\}).  
\end{equation}

\begin{proposition} \label{prop:52} Let $K\subset V$ be compact.  Then
$$
d_{\calE}(K) = d_{\calE}(\pi^{-1}(\pi(K))) = d(\pi(K))^{2},
$$
where $d(\cdot)$ denotes the classical transfinite diameter in $\CC^2$.
\end{proposition}

\begin{proof}  
The first equality follows immediately by applying a standard limiting argument to (\ref{eqn:vdme}).  For the second, we view $\calE$ as the standard monomial basis for $\CC[z_1,z_2]$, and hence $VDM(\cdot)$ is the same Vandermonde determinant that gives the classical transfinite diameter.  The exponent of 2 comes from taking a $(kM_k)$-th (rather than an $L_k$-th) root in the limiting process (see Proposition \ref{prop:46} with $n=2$).  
\end{proof}

\subsection{Relations between Chebyshev constants.}
Consider the complexified sphere $\bV(z_1^2+z_2^2+z_3^2-1)$ given projectively by $Z_1^2+Z_2^2+Z_3^2=Z_0^2$.  Consider 
local coordinates $v=(v_0,v_1)$ about the point  $$[Z_0:Z_1:Z_2:Z_3]=[v_0:v_1:1:v_3]=[0:0:1:i]\in\PP^3,$$
and let $\Delta$ be the Okounkov body calculated in these coordinates.  In these coordinates, $V=\bV(v_1^2+1+v_3^2-v_0^2)$, and one can calculate that $\Delta$ is again the triangle given by Proposition \ref{prop:oksphere}.  As a result of Theorem \ref{thm:Ok}, the integral formula 
\begin{equation} \label{eqn:Ok3}
d^w(K) = \exp\left( \frac{1}{\vol(\Delta)}\int_{\Delta^{\circ}} \log T^w(K,\theta)d\theta   \right)
\end{equation}
holds, where $w$ is a weight on a compact subset $K$ of $\{Z_0,Z_2\neq 0\}$.
We will show how the integrands in (\ref{eqn:thm51}) and (\ref{eqn:Ok3}) are related in a particular case.

First, we will transform (\ref{eqn:thm51}) to $v$-coordinates. Similar to the previous section, we use 
$v_1=z_1/z_2,\ v_0=1/z_2, \ v_3= z_3/z_2$, and polynomial evaluation is related by $v_0^{k}p(v_0,v_1) =p(z_1,z_2)$.  

Let us use tilded quantities ($\tilde\calE$, $\tilde\theta$, etc.) to denote the $v$-coordinate versions of the quantities in (\ref{eqn:thm51}), and use the same $\tilde\theta$ variable in equation (\ref{eqn:Ok3}); then $\tilde S\subset\Delta$.  With the weight $w(v)=v_0$, (\ref{eqn:thm51}) and (\ref{eqn:Ok3}) become  
\begin{eqnarray}
d_{\tilde\calE}^w(K) &=& \exp\left( \frac{1}{\vol(\tilde S)}\int_{\tilde S^{\circ}} \log T_{\tilde\calE}^w(K,\tilde\theta)\,d\tilde\theta  \right),    \label{eqn:dk1} \\
\quad d^w(K) &=& \exp\left(  \frac{1}{\vol(\Delta)} \int_{\Delta^{\circ}} \log T^w(K,\tilde\theta)\,d\tilde\theta  \right).  \label{eqn:dk2}
\end{eqnarray}
Figure 3 shows the relation between the different convex bodies and parameters. 

\begin{figure} \label{fig1}
\includegraphics[height=10cm]{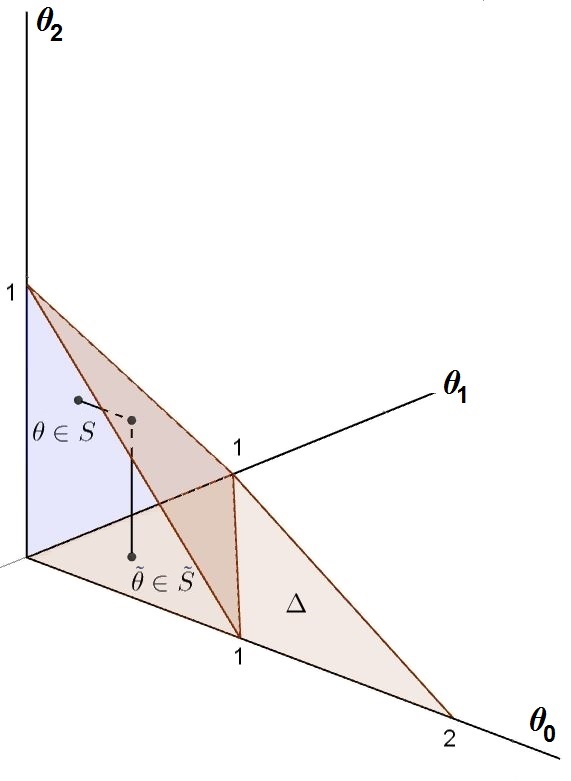}
\caption{Relations between the convex bodies. $S$ is constructed from monomials in affine coordinates, and $\tilde S\subset\Delta$ from monomials in $v$-coordinates.  The parameters are related by projection: we have $\theta=(\theta_1,\theta_2)$ and $\tilde\theta=(\theta_0,\theta_1)$ where $(\theta_0,\theta_1,\theta_2)$ is a point in the simplex $\theta_0+\theta_1+\theta_2=1$, $\theta_j>0\ \forall j$. } 
\end{figure}


\begin{definition}\rm
Let $A(r_1,r_2)$ be the annulus $\{\zeta\in\CC: r_1<|\zeta|<r_2\}$, with $r_1<1<r_2$, and let $U$ be an open subset of $V$. Suppose there exists a holomorphic map 
$$
A(r_1,r_2)\times U\ni (\zeta,v)\mapsto \zeta*v \in V
$$
such that $\zeta*(\eta*v) = (\zeta\eta)*v$ whenever $v,\eta*v\in U$ and $\zeta,\eta,\zeta\eta\in A(r_1,r_2)$.  The map $*$ is called a \emph{local circle action on $U$}, and $K\subset U$ is \emph{locally circled under $*$} if it is invariant under the restriction of the action to the unit circle:
$$
e^{i\theta}*v\in K \hbox{ if and only if } v\in K, \hbox{ for all }\theta\in\RR  \quad (\hbox{i.e. } e^{i\theta}*K=K).
$$
\end{definition}

A local circle action on a hypersurface $V\subset\CC^3$ arises naturally as follows. Locally (say in an open set $U$), $V$ is a graph over the $v_0,v_1$ variables; let us write
$v_3=\upsilon_3(v_0,v_1)$  for some holomorphic function $\upsilon_3$. Scalar multiplication $\zeta(v_0,v_1)=(\zeta v_0,\zeta v_1)$ lifts to a map
\begin{equation} \label{eqn:55}
\zeta*(v_1,v_2,v_3):=(\zeta v_1,\zeta v_2,\upsilon_3(\zeta v_1,\zeta v_2))
\end{equation}
which satisfies the above definition, as long as $\upsilon_3$ extends to a well-defined holomorphic function on some neighbourhood of a set of the form $\{e^{i\theta}w: \theta\in\RR\}$, for some $w\in\pi(U)$ (here $\pi$ denotes projection onto the $(v_0,v_1)$ variables).   This holds, for example, if $\upsilon_3$ has a Laurent series expansion at $w$, $\upsilon_3(w_0,w_1)=\sum_{i,j=-\infty}^{\infty}c_{ij}w_0^iw_1^j$. 
 
\medskip
 
 In particular, when $V$ is the sphere, let $r<1$ and consider $$W_r=\{(v_0,v_1)\in\CC^2: |v_0|,|v_1|\leq r\}.$$ In $v$-coordinates, we have $\upsilon_3(v_0,v_1)=(v_0^2-v_1^2-1)^{1/2}$.  
Then $\zeta*v$ given by (\ref{eqn:55}) is a local circle action defined on $$\tilde W_r = \{(v_0,v_1,\upsilon_3(v_0,v_1)): \  (v_0,v_1)\in W_r\}.$$ 
 
 
Let us see how a locally circled set can be generated in affine coordinates. 
 
 \begin{lemma} \label{lem:55} 
Let $K\subset\tilde W_r\subset V$ be a compact, locally circled set under the action $*$ on the sphere.  Using $*$, one can generate $K$ from a smaller set $K'\subset\CC\times[0,\infty)$ as follows:
$$
K = \{e^{i\theta}*(z_1,z_2,z_3)\in \tilde W_r: \  z_1=w, z_2=r, (w,r)\in K'\}.  
$$
\end{lemma}

\begin{proof}
Fix $v\in K$, which we will write as $z\in K$ in affine coordinates. We compute the holomorphic function $z(\zeta) = (z_1(\zeta),z_2(\zeta),z_3(\zeta))$ given by the transformation $\zeta\mapsto\zeta*v$:
$$z_1(\zeta) = \frac{\zeta v_1}{\zeta v_0} = \frac{v_1}{v_0} = z_1 ,\quad  z_2(\zeta) = \frac{1}{\zeta v_0} = \tfrac{1}{\zeta}z_2,
$$
and this determines $z_3(\zeta)$ by lifting to $\tilde W_r$.  Since $e^{i\theta}*v\in K$ for all $\theta$, we can choose $\theta_v$ such that $z_2(e^{i\theta_v}) = e^{-i\theta_v}z_2\in[0,\infty)$.  

We now vary $v$ and put $I=\{r=z_2(e^{i\theta_v}): v\in K\}\subset[0,\infty)$.  Then 
$$
K':=\{(w,r)\in\CC\times I:  (w,r,z_3)\in K \hbox{ for some } z_3\in\CC \}
$$
is the desired set.  
\end{proof}
This is a rotation about the origin in the $z_2$-plane, lifted to the variety.

\medskip

We return back to $v$-coordinates to relate Chebyshev constants. The notation in the following proposition is as in (\ref{eqn:dk1}) and (\ref{eqn:dk2}) above, with the weight given by $w(v_0,v_1)=v_0$.
 
\begin{proposition} \label{prop:55}
Let $r\in(0,1)$, and suppose $K\subset\tilde W_r\subset V$ is locally circled under $*$.  Then 
$$
T_{\tilde\calE}^w(K,\tilde\theta) = T^w(K,\tilde\theta)
$$
for all $\tilde\theta\in\tilde S\subset\Delta$. 
\end{proposition}

\begin{proof}
By definition, it is sufficient to show that for any $\alpha\in\NN$ with $|\alpha|\leq k$,
$$
\inf\{ \|p\|_{K,w}: p\in\CC[v_0,v_1],\,  \ttt(p)=v^{\alpha} \}  \ = \  \inf\{ \|p\|_{K,w}: p\in\CC[V], \,  \ttt(p)=v^{\alpha} \},
$$
or more compactly, $T^w_{\calE,k}(K,\alpha) = T^w_k(K,\alpha)$.  Note that any $p\in\CC[v_0,v_1]_{\leq k}$ is also a polynomial in $\CC[V]_{\leq k}$ with no term involving $v_3$, so the inf on the right-hand side is over a larger collection.  Hence $T^w_{\calE,k}(K,\alpha) \geq T^w_k(K,\alpha)$.

To prove the reverse inequality, fix $k\in\NN$ and $\alpha\in\NN_0^2$ with $|\alpha|\leq k$ (i.e., $\alpha/k\in\tilde S$).  Let $p\in\CC[V]_{\leq k}$ with $\ttt(p)=v^{\alpha}$, and for the moment, let $v=(v_0,v_1)$ be fixed in $K$.  Define  $\psi(\zeta):=p(\zeta*v)_{w,k}$; then $\psi$ is holomorphic on a neighborhood of the unit circle, as can be seen by writing it out in a series expansion:
\begin{eqnarray*}
\psi(\zeta) = p(\zeta*v)_{w,k} &=& \zeta^kv_0^k\sum_{j=0}^{\infty}\sum_{|\beta|=j} c_{\beta}(\zeta v_0)^{\beta_0}(\zeta v_1)^{\beta_1}  \\
&=& \sum_{j=0}^{\infty} \biggl(\sum_{|\beta|=j}c_{\beta}v_0^{\beta_0+k}v_1^{\beta_1} \biggr)\zeta^{j+k}.  
\end{eqnarray*}
Since $K$ is locally circled, $\zeta*v\in K$ for all $|\zeta|=1$, so that $|p(\zeta*v)|\leq\|p\|_K$.  Plugging this into the Cauchy integral formula for the  coefficient of $\zeta^{j+k}$, we have 
\begin{equation}\label{eqn:58}
\left|\sum_{|\beta|=j}c_{\beta}v_0^{\beta_0+k}v_1^{\beta_1}\right|  = \left| \frac{1}{2\pi i}\int_{|\zeta|=1} \frac{p(\zeta*v)}{\zeta^{j+k+1}}\, d\zeta \right| \leq \|p\|_K, \quad \forall\, j\in\NN.
\end{equation}
In particular, this is true when $j=|\alpha|\leq k$.  Since $v$ was an arbitrary point of $K$, let us now treat it as a variable and define the polynomial
$$\tilde p(v):=\sum_{|\beta|=|\alpha|}c_{\beta}v_0^{\beta_0}v_1^{\beta_1} \in \CC[v_1,v_2]_{\leq k}. $$
Clearly, by construction $\ttt(\tilde p)=v^{\alpha}$, and by (\ref{eqn:58}), $\|\tilde p\|_{K,w,k} \leq \|p\|_{K,w,k}$.   Hence
$$
T_{\calE,k}^w(K,\alpha) \leq \|\tilde p\|_{K,w,k} \leq \|p\|_{K,w,k}  
$$
and since $p\in\CC[V]_{\leq k}$ was an arbitrary polynomial with $\ttt(p)=v^{\alpha}$, we can take the inf over all such polynomials to obtain $T^w_{\calE,k}(K,\alpha) \leq T^w_k(K,\alpha)$.
\end{proof}

\begin{remark}\rm
The notion of a locally circled set is adapted from the notion of a \emph{circled set in $\CC^n$}.  
   Recall that a compact set $K\subset\CC^n$ is circled if $e^{i\theta}z\in K$ whenever $z\in K$.  For such sets, the {Chebyshev constants} and {homogeneous Chebyshev constants} (of Theorems \ref{thm:z} and \ref{thm:j} respectively) are equal. The proof is essentially the same as that of the above proposition.
\end{remark}



\bibliographystyle{abbrv}
\bibliography{myreferences}

\end{document}